\documentclass[twoside,11pt]{article}
\title{} \author{} \date{}
\usepackage{latexsym,amssymb,times}%,showkeys}
\input amssym.def
\newtheorem{te}{Theorem}[section]
\newtheorem{prop}[te]{Proposition}
\newtheorem{cor}[te]{Corollary}
\newtheorem{fac}[te]{Fact}
\newtheorem{lem}[te]{Lemma}
\newtheorem{rem}[te]{Remark}
\newtheorem{ex}[te]{Example}

\def\dok{\noindent{\bf Proof. }}
\def\kdok{\hfill $\Box$ \par \vspace*{2mm} }
\def\a{\alpha}
\def\b{\beta}
\def\g{\gamma}
\def\f{\varphi}

\def\o{\omega}

\def\s{\sigma}

\def\S{{\mathbb S}}
\def\P{{\mathbb P}}
\def\Q{{\mathbb Q}}

\def\N{{\mathbb N}}

\def\Z{{\mathbb Z}}

\def\H{{\mathcal H}}

\def\CS{{\mathcal S}}

\def\I{{\mathcal I}}

\def\down{\!\downarrow}

\def\la{\langle}
\def\ra{\rangle}

\def\Lim{\mathop{\mbox{Lim}}\nolimits}

\def\Fin{\mathop{\rm Fin}\nolimits}

\def\sm{\mathop{\rm sm}\nolimits}
\def\sq{\mathop{\rm sq}\nolimits}
\def\rp{\mathop{\rm rp}\nolimits}
\def\Block{\mathop{\rm Block}\nolimits}

\begin{document}
\thispagestyle{plain}
\begin{center}
           {\large \bf {\uppercase{Posets of copies of countable scattered\\[1mm] linear orders}}}
\end{center}
%\vspace{0.5mm}
\begin{center}
{\small \bf Milo\v s S.\ Kurili\'c}\\[1mm]
         {\small Department of Mathematics and Informatics, University of Novi Sad, \\
         Trg Dositeja Obradovi\'ca 4, 21000 Novi Sad, Serbia.\\[-1mm]
                                     e-mail: milos@dmi.uns.ac.rs}
\end{center}
%\vspace{0.5mm}
\begin{abstract}
\noindent
We show that the separative quotient of the poset $\la \P (L),\subset \ra$ of isomorphic suborders of a countable
scattered linear order $L$ is $\sigma$-closed and atomless. So, under the CH, all these posets are forcing-equivalent
(to $(P(\o )/\Fin )^+$).\footnote{
{\sl 2010 MSC}:
06A05, % Total order
06A06,  % Partial order, general
03C15,  % Denumerable structures
03E40, % Other aspects of forcing and Boolean valued models
03E35. % Consistency and independence results
\\
{\sl Key words and phrases}: scattered linear order, isomorphic substructure, denumerable structure, $\sigma$-closed poset, forcing.}
\end{abstract}

\section{Introduction}
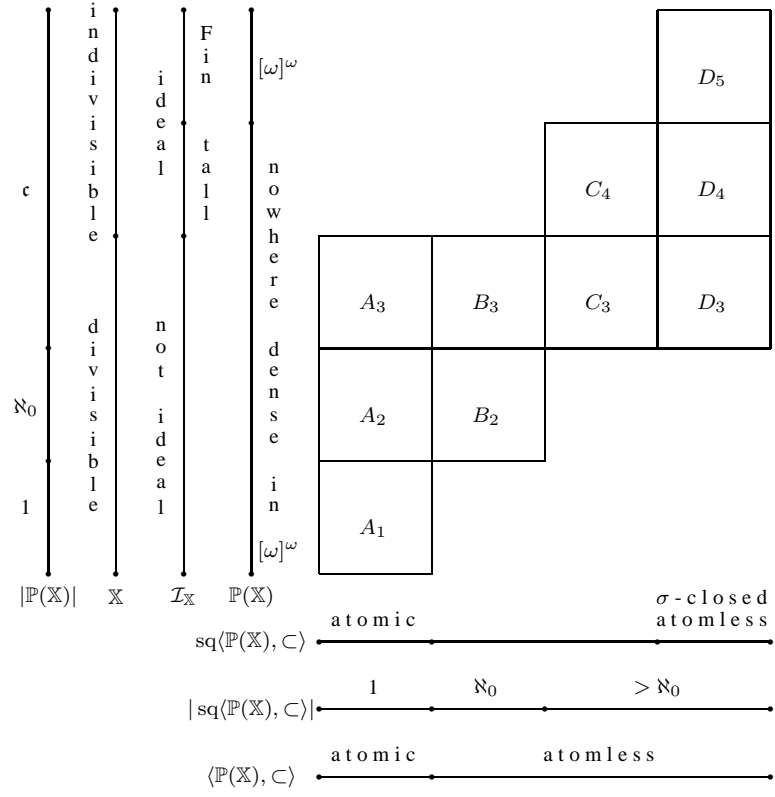
\begin{figure}[htb]
\begin{center}
\unitlength 0.6mm %0.7mm%1mm % = .854pt
\linethickness{0.4pt}
\ifx\plotpoint\undefined\newsavebox{\plotpoint}\fi % GNUPLOT compatibility

%================================  Perfect =============================

\begin{picture}(180,190)(0,0)

%----------------------------- linije --------------------------

\put(70,10){\line(1,0){100}}%1
\put(70,25){\line(1,0){100}}%2
\put(70,40){\line(1,0){100}}%3
\put(10,55){\line(0,1){125}}%4
\put(25,55){\line(0,1){125}}%5
\put(40,55){\line(0,1){125}}%6
\put(55,55){\line(0,1){125}}%7
\put(70,55){\line(0,1){75}}%8
\put(95,55){\line(0,1){75}}%9
\put(120,80){\line(0,1){50}}%10
\put(120,130){\line(0,1){25}}%11
\put(145,105){\line(0,1){75}}%12
\put(170,105){\line(0,1){75}}%13
\put(70,55){\line(1,0){25}}%14
\put(70,80){\line(1,0){50}}%15
\put(70,105){\line(1,0){100}}%16
\put(70,130){\line(1,0){100}}%17
\put(120,155){\line(1,0){50}}%18
\put(145,180){\line(1,0){25}}%19
%\put(70,200){\line(1,0){100}}%20

%----------------------------- tacke --------------------------

\put(70,10){\circle*{1}}%1
\put(95,10){\circle*{1}}%2
\put(170,10){\circle*{1}}%3
\put(70,25){\circle*{1}}%4
\put(95,25){\circle*{1}}%5
\put(120,25){\circle*{1}}%6
\put(170,25){\circle*{1}}%7
\put(70,40){\circle*{1}}%8
\put(145,40){\circle*{1}}%9
\put(170,40){\circle*{1}}%10
\put(10,55){\circle*{1}}%11
\put(10,80){\circle*{1}}%12
\put(10,105){\circle*{1}}%13
\put(10,180){\circle*{1}}%14
\put(25,55){\circle*{1}}%15
\put(25,130){\circle*{1}}%16
\put(25,180){\circle*{1}}%17
\put(40,55){\circle*{1}}%18
\put(40,130){\circle*{1}}%19
\put(40,155){\circle*{1}}%20
\put(40,180){\circle*{1}}%21
\put(55,55){\circle*{1}}%22
\put(55,155){\circle*{1}}%23
\put(55,180){\circle*{1}}%24
%\put(70,200){\circle*{1}}%25
%\put(95,200){\circle*{1}}%26
%\put(120,200){\circle*{1}}%27
%\put(145,200){\circle*{1}}%28
%\put(170,200){\circle*{1}}%29
\put(95,40){\circle*{1}}%1

%----------------------------- tekst -----------------------
%\small
%\footnotesize
\scriptsize
%\tiny
\put(82,15){\makebox(0,0)[cc]{a t o m i c}}%
\put(132,15){\makebox(0,0)[cc]{a t o m l e s s}}%
\put(82,45){\makebox(0,0)[cc]{a t o m i c}}%

\put(82,30){\makebox(0,0)[cc]{1}}%
\put(107,30){\makebox(0,0)[cc]{$\aleph _0$}}%
\put(145,30){\makebox(0,0)[cc]{$>\aleph _0$}}%

\put(157,50){\makebox(0,0)[cc]{$\sigma$ - c l o s e d}}%
\put(157,45){\makebox(0,0)[cc]{a t o m l e s s}}%

\put(5,70){\makebox(0,0)[cc]{1}}%
\put(5,92){\makebox(0,0)[cc]{$\aleph _0$}}%
\put(5,140){\makebox(0,0)[cc]{${\mathfrak c}$}}%
\put(20,180){\makebox(0,0)[cc]{i}}%
\put(20,175){\makebox(0,0)[cc]{n}}%
\put(20,170){\makebox(0,0)[cc]{d}}%
\put(20,165){\makebox(0,0)[cc]{i}}%
\put(20,160){\makebox(0,0)[cc]{v}}%
\put(20,155){\makebox(0,0)[cc]{i}}%
\put(20,150){\makebox(0,0)[cc]{s}}%
\put(20,145){\makebox(0,0)[cc]{i}}%
\put(20,140){\makebox(0,0)[cc]{b}}%
\put(20,135){\makebox(0,0)[cc]{l}}%
\put(20,130){\makebox(0,0)[cc]{e}}%
\put(20,110){\makebox(0,0)[cc]{d}}%
\put(20,105){\makebox(0,0)[cc]{i}}%
\put(20,100){\makebox(0,0)[cc]{v}}%
\put(20,95){\makebox(0,0)[cc]{i}}%
\put(20,90){\makebox(0,0)[cc]{s}}%
\put(20,85){\makebox(0,0)[cc]{i}}%
\put(20,80){\makebox(0,0)[cc]{b}}%
\put(20,75){\makebox(0,0)[cc]{l}}%
\put(20,70){\makebox(0,0)[cc]{e}}%
\put(35,165){\makebox(0,0)[cc]{i}}%
\put(35,160){\makebox(0,0)[cc]{d}}%
\put(35,155){\makebox(0,0)[cc]{e}}%
\put(35,150){\makebox(0,0)[cc]{a}}%
\put(35,145){\makebox(0,0)[cc]{l}}%
\put(35,110){\makebox(0,0)[cc]{n}}%
\put(35,105){\makebox(0,0)[cc]{o}}%
\put(35,100){\makebox(0,0)[cc]{t}}%
\put(35,90){\makebox(0,0)[cc]{i}}%
\put(35,85){\makebox(0,0)[cc]{d}}%
\put(35,80){\makebox(0,0)[cc]{e}}%
\put(35,75){\makebox(0,0)[cc]{a}}%
\put(35,70){\makebox(0,0)[cc]{l}}%
\put(45,175){\makebox(0,0)[cc]{F}}%
\put(45,170){\makebox(0,0)[cc]{i}}%
\put(45,165){\makebox(0,0)[cc]{n}}%
\put(45,150){\makebox(0,0)[cc]{t}}%
\put(45,145){\makebox(0,0)[cc]{a}}%
\put(45,140){\makebox(0,0)[cc]{l}}%
\put(45,135){\makebox(0,0)[cc]{l}}%
\put(61,167){\makebox(0,0)[cc]{$[\omega ]^{\omega }$}}%
\put(60,145){\makebox(0,0)[cc]{n}}%
\put(60,140){\makebox(0,0)[cc]{o}}%
\put(60,135){\makebox(0,0)[cc]{w}}%
\put(60,130){\makebox(0,0)[cc]{h}}%
\put(60,125){\makebox(0,0)[cc]{e}}%
\put(60,120){\makebox(0,0)[cc]{r}}%
\put(60,115){\makebox(0,0)[cc]{e}}%
\put(60,105){\makebox(0,0)[cc]{d}}%
\put(60,100){\makebox(0,0)[cc]{e}}%
\put(60,95){\makebox(0,0)[cc]{n}}%
\put(60,90){\makebox(0,0)[cc]{s}}%
\put(60,85){\makebox(0,0)[cc]{e}}%
\put(60,75){\makebox(0,0)[cc]{i}}%
\put(60,70){\makebox(0,0)[cc]{n}}%
\put(61,60){\makebox(0,0)[cc]{$[\omega ]^{\omega }$}}%

%\put(40,205){\makebox(0,0)[cc]{$\langle {\mathbb P} ({\mathbb X} ) , \subset \rangle$ i s}}%

%\put(82,205){\makebox(0,0)[cc]{1}}%
%\put(107,205){\makebox(0,0)[cc]{C o h e n}}%
%\put(157,195){\makebox(0,0)[cc]{u n d e r  (CH)}}%
%\put(157,205){\makebox(0,0)[cc]{$P(\o)/\Fin$}}%
%\put(40,198){\makebox(0,0)[cc]{f o r c i n g - e q u i v a l e n t }}%

\put(10,50){\makebox(0,0)[cc]{$|{\mathbb P} ({\mathbb X} )|$}}%
\put(25,50){\makebox(0,0)[cc]{${\mathbb X}$}}%
\put(40,50){\makebox(0,0)[cc]{${\mathcal I}_{\mathbb X}$}}%
\put(55,50){\makebox(0,0)[cc]{${\mathbb P}({\mathbb X})$}}%
\put(55,40){\makebox(0,0)[cc]{$\mathop{\rm sq}\nolimits  \langle {\mathbb P} ({\mathbb X} ) , \subset \rangle$}}%
\put(55,25){\makebox(0,0)[cc]{$|\mathop{\rm sq}\nolimits  \langle {\mathbb P} ({\mathbb X} ) , \subset \rangle|$}}%
\put(55,10){\makebox(0,0)[cc]{$\langle {\mathbb P} ({\mathbb X} ) , \subset \rangle$}}%

%--------------------A1D5------------------------------------------------

\put(82,65){\makebox(0,0)[cc]{$A_1$}}%
\put(82,90){\makebox(0,0)[cc]{$A_2$}}%
\put(82,115){\makebox(0,0)[cc]{$A_3$}}%
\put(107,90){\makebox(0,0)[cc]{$B_2$}}%
\put(107,115){\makebox(0,0)[cc]{$B_3$}}%
\put(132,115){\makebox(0,0)[cc]{$C_3$}}%
\put(132,140){\makebox(0,0)[cc]{$C_4$}}%
\put(157,115){\makebox(0,0)[cc]{$D_3$}}%
\put(157,140){\makebox(0,0)[cc]{$D_4$}}%
\put(157,165){\makebox(0,0)[cc]{$D_5$}}%
%\put(,){\makebox(0,0)[cc]{}}%
\end{picture}

\end{center}

\vspace{-7mm}

\caption{Binary relations on countable sets}\label{F4001}
\end{figure}

The posets of the form $\langle  {\mathbb P} ({\mathbb X} ), \subset\rangle $, where ${\mathbb X} $ is a relational structure and
${\mathbb P} ({\mathbb X} )$ the set of the domains of its isomorphic substructures, were investigated in \cite{Ktow}.
In particular, a classification of countable binary structures
related to the order-theoretic and forcing-related properties of the posets of their copies is described in Diagram \ref{F4001}:
for the structures from column $A$ (resp.\ $B$; $D$) the corresponding posets are forcing equivalent to the trivial poset
(resp.\ the Cohen forcing, $\langle  {}^{<\omega }2, \supset\rangle $;
an $\omega _1 $-closed atomless poset) and, for the structures from the class $C_4$, the posets of copies are forcing equivalent to the posets
of the form $(P(\omega )/{\mathcal I})^+$, for some co-analytic tall ideal ${\mathcal I}$.
For example, all countable non-scattered linear orders are in the class $C_4$, moreover, as a consequence of the main result of  \cite{KurTod} we have
\begin{te} \rm \label{T4120}
For each countable non-scattered linear order $L$ the poset $\la \P (L), \subset \ra$ is forcing equivalent to
the two-step iteration $\S \ast \pi$, where $\S$ is the Sacks forcing and $1_\S \Vdash `` \pi $ is a $\sigma$-closed forcing".
If the equality sh$(\S )=\aleph _1$
or PFA holds in the ground model, then the second iterand is forcing equivalent to the poset $(P(\o )/\Fin )^+$ of the
Sacks extension.
\end{te}
The aim of this paper is to complete the picture of countable linear orders in this context and,
having in mind Theorem \ref{T4120}, we concentrate our attention on countable scattered linear orders.
In the simplest case, if $L$ is the ordinal $\o$, then $\la \P (L), \subset \ra = \la [\o ]^\o ,\subset \ra$
is a homogeneous atomless partial order of size ${\mathfrak c}$ and its separative quotient,
the poset $(P(\o )/\Fin )^+$, is $\sigma$-closed. We will show that the same holds for
each countable scattered linear order. So the following theorem is our main result.
\begin{te}  \rm \label{T4121}
For each countable scattered linear order $L$ the poset
$\la \P (L), \subset \ra$ is homogeneous, atomless, of size ${\mathfrak c}$ and its separative quotient is
$\sigma$-closed.
\end{te}
\begin{cor} \rm \label{T4118}
If $L$ is a countable linear order,
then the poset $\la \P (L), \subset \ra$ is forcing equivalent to

- $\S \ast \pi$, where
  $1_\S \Vdash `` \pi $ is $\sigma$-closed", if $L$ is non-scattered \cite{KurTod};

- A $\sigma$-closed atomless forcing, if $L$ is scattered.

\noindent
Under the CH,  the poset $\la \P (L), \subset \ra$ is forcing equivalent to

- $\S \ast \pi $, where $1_\S \Vdash `` \pi =(P(\check{\o })/\Fin )^+"$, if $L$ is non-scattered \cite{KurTod};

- $(P(\o )/\Fin )^+$, if $L$ is scattered.
\end{cor}
The most difficult part of the proof of Theorem \ref{T4121} is to show that the separative quotient of $\la \P (L), \subset \ra$ is $\sigma$-closed
(this result is the best possible: ``$\s$-closed" can not be replaced by ``$\o _2$-closed", see Example \ref{EX4001}).
Namely, it is easy to see that there are copies of an $\o$-sum $\sum _\o L_i$ of linear orders $L_i$, which are not of the form
$\bigcup _{i\in \o}C_i$, where $C_i \in \P (L_i)$, so the Hausdorff hierarchy of scattered linear orders
can not  be used (easily) for an inductive proof.
Instead of that hierarchy we use the result of Laver \cite{Lav} that a countable scattered linear order is a finite sum of
hereditarily indecomposable (ha) linear orders. So we first prove the statement for ha-orders, then for special blocks of
ha-orders and, finally, for finite sums of blocks.
\section{Preliminaries}
A linear order $L$ is said to be {\it scattered} iff it does not contain a dense suborder or, equivalently,
iff the rational line, $\Q$, does not embed in $L$.
By $\CS$ we denote the class of all countable scattered linear orders.
\begin{fac}  \rm \label{T4102}
If $L$ is a linear order satisfying $L+L\hookrightarrow L$, then $L$ is not scattered (see \cite{Rosen}, p.\ 180).
\end{fac}
\dok
By the assumption,  $L+(L+L)\hookrightarrow L+L\hookrightarrow L$. By recursion we construct the sequences
$\la L_{\f }: \f \in {}^{<\o }2\ra$ and $\la L_{\f }': \f \in {}^{<\o }2\ra$ in $\P (L)$ and
$\la q_{\f }: \f \in {}^{<\o }2\ra$ in $L$ such that
(i) $L_{\emptyset}=L$,
(ii) $L_{\f ^\smallfrown 0} < L_\f '< L_{\f ^\smallfrown 1}$,
(iii) $L_{\f ^\smallfrown 0} \cup  L_\f '\cup L_{\f ^\smallfrown 1} \subset L_\f$,
(iv) $q_\f \in L_\f '$.
Then $\{ q_{\f }: \f \in {}^{<\o }2\}$ is a copy of $\Q$ in $L$.
\kdok
A linear order $L$ is said to be {\it additively indecomposable}
(respectively {\it left indecomposable; right indecomposable}) iff for each decomposition
$L=L_0 + L_1$ we have $L\hookrightarrow L_0$ or $L\hookrightarrow L_1$
(respectively $L\hookrightarrow L_0$; $L\hookrightarrow L_1$). The class $\H$ of {\it hereditarily additively
indecomposable} (or {\it ha-indecomposable}) linear orders is the smallest class of order types of
countable linear orders containing the one element order type, {\bf 1}, and containing the $\o$-sum, $\sum _\o L_i$,
and the $\o ^*$-sum, $\sum _{\o ^*} L_i$, for each sequence $\la L_i :i\in \o \ra$ in $\H$ satisfying
\begin{equation}\label{EQ4100}
\forall i\in \o \;\; |\{ j\in \o : L_i \hookrightarrow L_j \}|=\aleph _0 .
\end{equation}
\begin{fac}  \rm \label{T4100}
(a) $\H \subset \CS$ (see \cite{Rosen}, p.\ 196);

(b) If $L\in \H$ is an $\o$-sum, then $L$ is right indecomposable (see \cite{Rosen}, p.\ 196);

(c) If $L\in \H$ is an $\o ^*$-sum, then $L$ is left indecomposable (see \cite{Rosen}, p.\ 196);

(d) If $L\in \CS$ is additively indecomposable, then $L$ is left indecomposable or right indecomposable
(see \cite{Rosen}, p.\ 175);

(e) (Laver, \cite{Lav}) If $L\in \CS$, then $L\in \H$ iff $L$ is additively indecomposable (see \cite{Rosen}, p.\ 201);

(f) (Laver, \cite{Lav}) If $L\in \CS$, then $L$ is a finite sum of elements of $\H$ (see \cite{Rosen}, p.\ 201).
\end{fac}
Let ${\mathbb P} =\langle  P , \leq \rangle $ be a pre-order. Then
$p\in P$ is an {\it atom} iff each  $q,r\leq p$ are compatible (there is $s\leq q,r$).
${\mathbb P} $ is called: {\it atomless} iff it has no atoms;
{\it homogeneous} iff  it has the largest element and $\P \cong p\down $, for each $p\in \P$.
If $\kappa $ is a regular cardinal, ${\mathbb P} $ is called $\kappa ${\it -closed} iff for each
$\gamma <\kappa $ each sequence $\langle  p_\alpha :\alpha <\gamma\rangle $ in $P$, such that $\alpha <\beta \Rightarrow p_{\beta}\leq p_\alpha $,
has a lower bound in $P$. $\omega _1$-closed pre-orders are called {\it $\sigma$-closed}.
Two pre-orders ${\mathbb P}$ and ${\mathbb Q}$ are called {\it forcing equivalent} iff they produce the same generic extensions.
\begin{fac}\rm\label{T4043}
If ${\mathbb P} _i$, $i\in I$, are $\kappa $-closed pre-orders, then $\prod _{i\in I}{\mathbb P} _i$ is $\kappa $-closed.
\end{fac}
A partial order ${\mathbb P} =\langle  P , \leq \rangle $ is called
{\it separative} iff for each $p,q\in P$ satisfying $p\not\leq q$ there is $r\leq p$ such that $r \perp q$.
The {\it separative modification} of ${\mathbb P}$
is the separative pre-order $\mathop{\rm sm}\nolimits ({\mathbb P} )=\langle  P , \leq ^*\rangle $, where
\begin{equation}\label{EQ4141}
p\leq ^* q \Leftrightarrow \forall r\leq p \; \exists s \leq r \; s\leq q .
\end{equation}
The {\it separative quotient} of ${\mathbb P}$
is the separative partial order $\mathop{\rm sq}\nolimits  ({\mathbb P} )=\langle P /\!\! =^* , \trianglelefteq \rangle$, where
$p = ^* q \Leftrightarrow p \leq ^* q \land q \leq ^* p\;$ and $\;[p] \trianglelefteq [q] \Leftrightarrow p \leq ^* q $.
\begin{fac}  \rm \label{T2226}
Let ${\mathbb P} , {\mathbb Q} $ and ${\mathbb P} _i$, $i\in I$, be partial orderings. Then

(a) ${\mathbb P}$, $\mathop{\rm sm}\nolimits ({\mathbb P})$ and $\mathop{\rm sq}\nolimits  ({\mathbb P})$ are forcing equivalent forcing notions;

(b) $\mathop{\rm sm}\nolimits ({\mathbb P} )$ is $\kappa $-closed iff $\mathop{\rm sq}\nolimits  ({\mathbb P} )$ is $\kappa $-closed;

(c) If $p_0, p_1, \dots p_n \in \P$, where $p_n \leq ^*  p_{n-1} \leq ^* \dots \leq ^* p_0$,
then there is $q\in \P$ such that $q \leq p_k$, for all $k\leq n$.

(d) ${\mathbb P} \cong {\mathbb Q}$ implies that $\mathop{\rm sm}\nolimits {\mathbb P} \cong  \mathop{\rm sm}\nolimits {\mathbb Q}$ and $\mathop{\rm sq}\nolimits  {\mathbb P} \cong  \mathop{\rm sq}\nolimits  {\mathbb Q}$;

(e) $\mathop{\rm sm}\nolimits (\prod _{i\in I}{\mathbb P} _i) = \prod _{i\in I}\mathop{\rm sm}\nolimits {\mathbb P} _i$ and
 $\mathop{\rm sq}\nolimits  (\prod _{i\in I}{\mathbb P} _i) \cong \prod _{i\in I}\mathop{\rm sq}\nolimits  {\mathbb P} _i$.

(f) If $X$ is an infinite set, $\I \subset P(X)$ an ideal containing $[X]^{<\o }$ and $\I ^+ = P(X)\setminus \I$ the corresponding family of $\I$-positive sets, then
$\sm \la \I ^+ , \subset \ra $ $= \la \I ^+ , \subset _\I \ra$, where $A\subset _\I B \Leftrightarrow A\setminus B \in \I$, for $A,B \in \I ^+$.
Also $\sq \la \I ^+ , \subset \ra =(P(X)/\I )^+$.
%(g?) ${\mathbb P}$ is atomless iff $\mathop{\rm sm}\nolimits ({\mathbb P} )$ is atomless iff $\mathop{\rm sq}\nolimits  ({\mathbb P} )$ is atomless;
\end{fac}
\dok
All the statements are folklore except, maybe, (c).
For a proof of (c), by recursion we define the sequence $\la q_k : k\leq m\ra$ such that
(i) $q_0=p_n$ and
(ii) $q_k \leq q_{k-1}, p_{n-k}$, for $0<k\leq n$.
Then $q_n\leq p_k$, for all $k\leq n$.
\kdok
\begin{fac}   \rm \label{T4056}
(Folklore) Under the CH, each  atomless separative $\omega _1$-closed pre-order of size $\omega _1$
is forcing equivalent to $(P(\omega )/\mathop{\rm Fin}\nolimits )^+$.
\end{fac}
We recall that the ideal $\mathop{\rm Fin}\nolimits  \times \mathop{\rm Fin}\nolimits  \subset P(\omega\times \omega )$
is defined by:
$$
\Fin \times \Fin =\{ A\subset \o \times \o : | \{ i\in \o : |A \cap L_i | =\o \} | < \o \},
$$
where $L_i=\{ i\} \times \o$, for $i\in \o$.
By ${\mathfrak h}({\mathbb P} )$ we denote the {\it distributivity number} of a poset ${\mathbb P}$. In particular, for $n\in {\mathbb N}$, let
${\mathfrak h}_n={\mathfrak h}(((P(\omega )/\mathop{\rm Fin}\nolimits  )^+ )^n )$; thus ${\mathfrak h}={\mathfrak h}_1$.
\begin{fac}    \rm  \label{T4091}
(a) $\mathop{\rm sm}\nolimits (\langle  [\omega ]^\omega , \subset \rangle ^n ) =\langle  [\omega ]^\omega , \subset ^*\rangle ^n$ and
$\mathop{\rm sq}\nolimits  (\langle  [\omega ]^\omega , \subset \rangle ^n) =((P(\omega )/\mathop{\rm Fin}\nolimits  )^+ )^n$ are forcing equivalent,
${\mathfrak t}$-closed atomless  pre-orders of size ${\mathfrak c}$.

(b) (Shelah and Spinas \cite{SheSpi1}) Con(${\mathfrak h}_{n+1}  < {\mathfrak h}_n$), for each $n\in {\mathbb N}$.

(c) (Szyma\'nski and Zhou \cite{Szym}) $(P(\omega \times \omega)/(\mathop{\rm Fin}\nolimits  \times \mathop{\rm Fin}\nolimits  ))^+$
is an $\omega _1$-closed, but not $\omega _2$-closed atomless  poset.

(d) (Hern\'andez-Hern\'andez \cite{Her})
Con(${\mathfrak h}((P(\omega \times \omega )/(\mathop{\rm Fin}\nolimits  \times \mathop{\rm Fin}\nolimits  ))^+) < {\mathfrak h}$).
\end{fac}
Now we prove the first part of Theorem \ref{T4121}.
\begin{prop}  \rm \label{T4117}
For each countable scattered linear order $L$ the partial ordering
$\la \P (L), \subset \ra$ is homogeneous, atomless and of size ${\mathfrak c}$.
\end{prop}
\dok
The homogeneity of $\la \P (L), \subset \ra$ is evident. For a proof that it is atomless first we show
\begin{equation}\textstyle \label{EQ4140}
\forall L\in \H \;\; (|L|=\o \Rightarrow \exists X,Y\in \P (L) \;\; X\cap Y=\emptyset ).
\end{equation}
If $L$ is an $\o$-sum, that is $L=\sum _{\o }L_i$, where
$\la L_i :i\in \o \ra$ is a sequence in $\H$ satisfying (\ref{EQ4100}), by recursion we define the sequences
$\la k_i :i\in \o \ra$ and $\la l_i :i\in \o \ra$ in $\o$ such that for each $i$

(i) $k_i <l_i$,

(ii) $l_i < k_{i+1}$,

(iii) $L_i \hookrightarrow L_{k_i}, L_{l_i}$.

\noindent
Using (\ref{EQ4100}) we choose $k_0, l_0 \in \o$ such that $k_0 <l_0$ and $L_0 \hookrightarrow L_{k_0}, L_{l_0}$.

Let the sequences $k_0, \dots , k_i$ and $l_0, \dots , l_i$ satisfy (i)-(iii). Then
$k_0 < l_0 < \dots <k_i <l_i$.
Using (\ref{EQ4100}) we choose $k_{i+1}, l_{i+1} \in \o$ such that $l_i <k_{i+1} <l_{i+1}$
and $L_{i+1} \hookrightarrow L_{k_{i+1}}, L_{l_{i+1}}$. Thus, the recursion works.

By (iii) there are $X_i,Y_i\cong L_i$ such that $X_i\subset L_{k_i}$ and $Y_i \subset L_{l_i}$.
Then $X=\sum _{\o }X_i$, $Y=\sum _{\o }Y_i \cong L$ and, by (i) and (ii) we have $X\cap Y=\emptyset$.

If $L$ is an $\o ^*$-sum, we proceed in the same way. Thus (\ref{EQ4140}) is proved.

By Fact \ref{T4100} for $L\in \CS$ there is $m\in \N$ such that $L=\sum _{i<m}L_i$, where $L_i\in \H$.
Let $J=\{ i<m : |L_i|=\o \}$. By (\ref{EQ4140}), for $i\in J$ there are $X_i,Y_i \in \P (L_i)$ such that
$X_i \cap Y_i =\emptyset$.
Let $X= \bigcup _{i\in J}X_i \cup \bigcup _{i\in m\setminus J}L_i$ and
$Y= \bigcup _{i\in J}Y_i \cup \bigcup _{i\in m\setminus J}L_i$.
Then $X,Y\in \P (L)$ and $|X\cap Y|= |\bigcup _{i\in m\setminus J}L_i|<\omega $ and, hence,
$X$ and $Y$ are incompatible elements of the poset $\la \P (L), \subset \ra$. So, since $\la \P (L), \subset \ra$
is a homogeneous partial order, it is atomless.

It is known (see \cite{Fra}, p.\ 170) that the equivalence classes corresponding to the relation $\sim$ on $L$,
defined by $x\sim y$ iff $|[\min \{ x,y \}, \max \{ x,y \}]|<\o$,
are convex parts of $L$ which are finite or isomorphic to $\o$, or $\o ^*$ or $\Z$.
Since $|L|=\o$ and two consecutive parts
can not be finite, there is one infinite part, say $L'$, and, clearly, it has ${\mathfrak c}$-many copies.
For each $C\in \P (L')$ we have $(L\setminus L')\cup C \in  \P(L)$ and, hence, $|\P (L)|={\mathfrak c}$.
\kdok
In the rest of the paper we prove that $\sq \la \P (L), \subset \ra$ is a $\s$-closed poset, for each countable scattered linear order $L$.
By Fact \ref{T2226}(b), it is sufficient to show that the pre-order $\sm \la \P (L), \subset \ra$ is $\s$-closed.
In the sequel we use the following notation:
$$\sm \langle \P (L) , \subset \rangle =\langle \P (L) , \leq \rangle.$$
\section{Elements of $\H$}
\begin{prop}\rm\label{T4103}
Let $L=\sum _\o L_i \in \H$, where $\la L_i :i\in \o \ra$ is a sequence in $\H$ satisfying (\ref{EQ4100}). Then

(a) $A\subset L$ contains a copy of $L$ iff
for each $ i, m\in \o $ there is finite $K \subset \o \setminus m$ such that
$L_i \hookrightarrow \bigcup _{j\in K }L_j \cap A $.
So, each $A\in \P (L)$ intersects infinitely many $L_i$'s.

(b) If $A,B\in \P (L)$, then $A\leq B$ iff for each $C\in \P (L)$ satisfying $C\subset A$
and each $i, m\in \o $ there exists a finite $K \subset \o \setminus m$ such that
$L_i \hookrightarrow \bigcup _{j\in K }L_j \cap C \cap B $.

(c) sm$\langle \P (L) , \subset \rangle$ is a $\sigma$-closed pre-order.

\noindent
The same statement holds for the $\o ^*$-sum $\sum _{\o ^*} L_i $.
\end{prop}

\dok
(a) ($\Rightarrow$) Let $f: L\hookrightarrow L$ and $C=f[L]\subset A$. Then $C=\sum _{i\in \o} f[L_i]$.

\vspace{2mm}
\noindent
{\it Claim 1}.
For each $i\in \o$ there is a finite set $K\subset \o$ such that $f[L_i]\subset \bigcup _{j\in K} L_j$.

\vspace{2mm}
\noindent
{\it Proof of Claim 1}. Since $f$ is an embedding and $L_i <L_{i+1}$ we have $f[L_i] <f[L_{i+1}]$.
For $x\in L_{i+1}$ we have $f(x) \in f[L_{i+1}]\subset \bigcup _{j\in \o} L_j$ and, hence, $f(x)\in L_{j_0}$, for some
$j_0 \in \o$. Now, by the monotonicity of $f$ we have $f[L_i]< \{ f(x) \}\subset L_{j_0}$, thus
$f[L_i]\subset \bigcup _{j\leq j_0}L_j$, so we can take $K=j_0 +1$ and Claim 1 is proved.

For $i\in \o$ let $K_i =\{ j\in \o : f[L_i]\cap L_j \neq \emptyset \}$. By Claim 1 we have
\begin{equation}\label{EQ4107}\textstyle
K_i \in [\o]^{<\o } \mbox{ and } f[L_i]\subset \bigcup _{j\in K_i}L_j.
\end{equation}

\noindent
{\it Claim 2}.
$K_i \leq K_{i+1}$, for each $i\in \o$. Consequently, either $K_i \cap K_{i+1}=\emptyset$ or $K_i \cap K_{i+1}=\{ \max K_i\}=\{ \min K_{i+1}\}$.

\vspace{2mm}
\noindent
{\it Proof of Claim 2}. Let $j'\in K_i$ and $j'' \in K_{i+1}$. Then there are $x\in L_i$ and $y\in L_{i+1}$ such that
$f(x)\in L_{j'}$ and $f(y)\in L_{j''}$ and, clearly, $x<y$. Now $j''<j'$ would imply $f(y)<f(x)$, which is impossible.
Thus $j'\leq j''$. Claim 2 is proved.

\vspace{2mm}
\noindent
{\it Claim 3}. $\bigcup _{i\in \o}K_i$ is an infinite subset of $\o$.

\vspace{2mm}
\noindent
{\it Proof of Claim 3}. On the contrary, suppose that $j_0=\max \bigcup _{i\in \o}K_i$. Let
$i_0 =\min \{ i\in \o : j_0 \in K_i \}$. Then $j_0 \in K_{i_0}\leq \{ j_0 \}$ and, by Claim 2,
$$
\forall i>i_0 \;\; (K_i=\{ j_0 \} \land f[L_i]\subset L_{j_0}) .
$$
By (\ref{EQ4100}), there are $i_1 , i_2 \in \o$ such that $i_0 +1 <i_1 <i_2$ and
$L_{j_0}\hookrightarrow L_{i_1}, L_{i_2}$, which implies
$L_{j_0}+L_{j_0}\hookrightarrow L_{i_1} + L_{i_2}\hookrightarrow f[L_{i_1}] + f[L_{i_2}] \subset L_{j_0}$.
But $L_{j_0}$ is a scattered linear order and, by Fact \ref{T4102}, $L_{j_0}+L_{j_0}\not\hookrightarrow L_{j_0}$.
A contradiction. Claim 3 is proved.

Let $i_0,m_0 \in \o$. By (\ref{EQ4100}), the set $I_{i_0}=\{ j\in \o : L_{i_0}\hookrightarrow L_j \}$ is an infinite set.

\vspace{2mm}
\noindent
{\it Claim 4}.
There is $j_0\in I_{i_0}$ such that $K_{j_0}\cap m_0=\emptyset$.

\vspace{2mm}
\noindent
{\it Proof of Claim 4}.
On the contrary, suppose that $K_j \cap m_0\neq\emptyset$, for each $j\in I_{i_0}$. Then
\begin{equation}\label{EQ4108}\textstyle
\forall j\in I_{i_0}\;\;\min K_j <m_0.
\end{equation}
For $i\in \o$ there is $j\in I_{i_0}$ such that $j>i+1$ and, by Claim 2, $K_i \leq K_{i+1}\leq K_j$ and, by (\ref{EQ4108}),
$\max K_i \leq \min K_{i+1}\leq \min K_j <m_0$. Thus $K_i\subset m_0$, for all $i\in \o$, which is impossible by Claim 3.
Claim 4 is proved.

By Claim 4, $K_{j_0}\in [\o \setminus m_0]^{<\o }$. By (\ref{EQ4107}) we have
$f[L_{j_0}]\subset \bigcup _{j\in K_{j_0}}L_j$. Since $j_0\in I_{i_0}$ and $f[L_{j_0}]\subset C\subset A$ we have
$L_{i_0}\hookrightarrow L_{j_0}\hookrightarrow f[L_{j_0}] \subset \bigcup _{j\in K_{j_0}}L_j \cap A$ and the proof of
``$\Rightarrow$" is finished.

($\Leftarrow$) Suppose that a set $A\subset L$ satisfies the given condition. By recursion we define
the sequences $\la K_i :i\in \o \ra $ and $\la f_i :i\in \o \ra $ such that for each $i\in \o$

(i) $K_i \in [\o ]^{<\o }$,

(ii) $K_0 < K_1 < \dots$,

(iii) $f_i : L_i \hookrightarrow \bigcup _{j\in K_i}L_j \cap A$.

\noindent
By the assumption, for $i=m=0$ there are $K_0 \in [\o ]^{<\o }$ and
$f_0 : L_0 \hookrightarrow \bigcup _{j\in K_0}L_j \cap A$.

Let $K_0, \dots , K_i$ and $f_0, \dots , f_i$ satisfy (i)-(iii) and let $m=\max (\bigcup _{r\leq i}K_r)+1$.
By the assumption for $i+1$ and $m$ there are $K_{i+1} \in [\o \setminus m ]^{<\o }$ and
$f_{i+1} : L_{i+1} \hookrightarrow \bigcup _{j\in K_{i+1}}L_j \cap A$ and the recursion works.

Let $f=\bigcup _{i\in \o }f_i$. By (ii) and (iii), $i_1 <i_2$ implies $K_{i_1}<K_{i_2}$, which implies
$f_{i_1}[L_{i_1}]<f_{i_2}[L_{i_2}]$ and, hence, $f:L\hookrightarrow A$. Thus $C=f[L]\in \P (L)$ and $C\subset A$.

(b) By (\ref{EQ4141}), $A\leq B$ iff for each $C\in \P (L)$ satisfying $C\subset A$
the set $C\cap B$ contains a copy of $L$. Now we apply (a) to $C\cap B$.

(c) For $A_n \in \P (L)$, $n\in \o$, where $A_0 \geq A_1 \geq \dots$ we will construct $A \in \P (L)$ such that
$A\leq A_n$, for all $n\in \o$.
First, by Fact \ref{T2226}(c), there are $C_i \in \P (L)$,
$i\in \o$, such that $C_0=A_0$ and
\begin{equation}\label{EQ4109}\textstyle
\forall i\in \o \;\; C_i \subset A_0 \cap \dots \cap A_i .
\end{equation}
By recursion we define
the sequences $\la K_i :i\in \o \ra $ and $\la f_i :i\in \o \ra $ such that for each $i\in \o$

(i) $K_i \in [\o ]^{<\o }$,

(ii) $K_i < K_{i+1} $,

(iii) $f_i : L_i \hookrightarrow \bigcup _{j\in K_i}L_j \cap C_i$.

\noindent
Since $C_0=A_0\in \P (L)$, by (a), for $i=m=0$ there are $K_0 \in [\o ]^{<\o }$ and
$f_0 : L_0 \hookrightarrow \bigcup _{j\in K_0}L_j \cap C_0$.

Let the sequences $K_0, \dots , K_{i'}$ and $f_0, \dots , f_{i'}$ satisfy (i)-(iii).
Since $A_{i'+1}\leq A_{i'}$, $C_{i'+1}\in \P (L)$
and, by (\ref{EQ4109}), $C_{i'+1}\subset A_{i'+1}$, according to (b), for $i'+1$ and
$m=\max (K_0 \cup \dots \cup K_{i'})+1$ there are
\begin{equation}\label{EQ4110}\textstyle
K_{i'+1} \in [\o \setminus (\max (K_0 \cup \dots \cup K_{i'})+1 )]^{<\o }
\end{equation}
\begin{equation}\label{EQ4111}\textstyle
f_{i'+1} : L_{i'+1} \hookrightarrow \bigcup _{j\in K_{i'+1}}L_j \cap C_{i'+1}
\end{equation}
(since, by (\ref{EQ4109})), $C_{i'+1} \cap A_{i'}=C_{i'+1}$). By (\ref{EQ4110})) we have (i) and (ii) and (iii) follows
from (\ref{EQ4111})). The recursion works.

Let $f=\bigcup _{i\in \o }f_i$. By (ii) and (iii), $i_1 <i_2$ implies $K_{i_1}<K_{i_2}$, which implies
$f_{i_1}[L_{i_1}]<f_{i_2}[L_{i_2}]$ and, hence, $f:L\hookrightarrow L$.
Thus
\begin{equation}\label{EQ4112}\textstyle
A=f[L]=\bigcup _{i\in \o }f_i[L_i]\in \P (L).
\end{equation}
Using the characterization from (b), for $n^* \in \o$ we show that $A\leq A_{n^*}$.
So, for $C^*\in \P (L)$ such that $C^* \subset A$ and $i^*,m^* \in \o$ we prove that
\begin{equation}\label{EQ4113}\textstyle
\exists K \in [\o \setminus m^* ]^{<\o }\;\;L_{i^*}\hookrightarrow \bigcup _{j\in K}L_j \cap C^* \cap A_{n^*}.
\end{equation}
By (ii), (iii) and (\ref{EQ4112}) we have $A=\sum _{i\in \o}\Lambda _i \cong L$, where $\Lambda _i =f_i[L_i]\cong L_i$,
thus $A\in \H$. Since $C^* \cong L \cong A$ we have $C^* \in \P (A)$ so, applying (a) to the linear order
$A$ instead of $L$ we obtain
\begin{equation}\label{EQ4114}\textstyle
\forall i,m \in \o \;\; \exists K \in [\o \setminus m ]^{<\o }\;\;f_i[L_i]\hookrightarrow \bigcup _{j\in K}f_j[L_j] \cap C^* .
\end{equation}
Let $m' >m^* , n^*$. By (\ref{EQ4114}), for $i^*$ and $m'$ there is
\begin{equation}\label{EQ4115}\textstyle
K^* \in [\o \setminus m']^{<\o } \mbox{ such that}
\end{equation}
\begin{equation}\label{EQ4116}\textstyle
f_{i^*}[L_{i^*}] \hookrightarrow \bigcup _{j\in K^*}f_j[L_j] \cap C^* .
\end{equation}
By (\ref{EQ4115}), for $j\in K^*$ we have $j>n^*$ and, by (\ref{EQ4109}), $C_j \subset A_{n^*}$. Thus, by (iii) we have
$f_j[L_j]\subset \bigcup _{s\in K_j}L_s \cap C_j \subset \bigcup _{s\in K_j}L_s \cap A_{n^*}$ which, together with
(iii) and (\ref{EQ4116}) gives
$L_{i^*} \hookrightarrow f_{i^*}[L_{i^*}]\hookrightarrow  \bigcup _{j\in K^*}f_j[L_j] \cap C^*
\subset  \bigcup _{j\in K^*} \bigcup _{s\in K_j}L_s \cap A_{n^*} \cap C^*
=\bigcup _{s\in \bigcup _{j\in K^*} K_j}L_s \cap C^* \cap A_{n^*}$.

In order to finish the proof of (\ref{EQ4113}) we prove that $\bigcup _{j\in K^*} K_j \cap m^* =\emptyset$.
By (\ref{EQ4115}), for $j\in K^*$ we have $j>m^*$. By (ii) the sequence $\la \min K_i : i\in \o \ra$ is increasing
and, hence, $\min K_j \geq j >m^*$, which implies $K_j \cap m^* =\emptyset$ and (\ref{EQ4113}) is proved.
\kdok

\section{Finite sums of $\o$-sums. Finite sums of $\o ^*$-sums}

\begin{lem}\rm\label{T4105}
Let $L_0=\sum _{\o } L_i^0 , L_1=\sum _{\o } L_i^1 \in \H$, where $\la L_i^0 :i\in \o \ra$  and $\la L_i^1 :i\in \o \ra$
are sequences in $\H$ satisfying (\ref{EQ4100}). Then

(a) $\exists i\in \o \;\; L_0 \hookrightarrow L_i ^1 \Leftrightarrow
\exists m \in \o \;\; L_0 \hookrightarrow \sum _{i\leq m}L_i ^1$;

(b) $L_0 +L_1 \not\in \H \Rightarrow \neg \exists m\in \o \;\;L_0 \hookrightarrow \sum _{i\leq m}L_i ^1$.

(c) If $L=L_0 +L_1 \not\in \H $ and $f:L\hookrightarrow L$, then $f[L_k]\subset L_k$, for $k=0,1$.
\end{lem}
\dok
(a) Suppose that $L_0 \hookrightarrow \sum _{i\leq m}L_i ^1$ and let $i_0=\max \{i\leq m : f[L_0]\cap L_i^1\neq\emptyset\}$.
Then $f[L_0]\cap L_{i_0}^1$ is a final part of the ordering $f[L_0]\cong L_0$ and, by Fact \ref{T4100}(a),
contains a copy of $L_0$. Thus $L_0\hookrightarrow L_{i_0}^1$.

(b) If $L_0 \hookrightarrow \sum _{i\leq m}L_i ^1$ then, by (a), there are $i_0\in \o$ and
$f: L_0\hookrightarrow L_{i_0}^1$. Then $\la L_0, L_0^1 , L_1^1 , \dots , L_{i_0}^1 ,\dots \ra$ is a sequence in $\H$
satisfying (\ref{EQ4100}) and $L_0 + L_1=L_0 + L_0^1 + L_1^1 + \dots + L_{i_0}^1 +\dots \in \H$.

(c) Suppose that $f[L_0]\cap L_1\neq\emptyset$. Then $f[L_0]\cap L_1$ is a final part of the ordering $f[L_0]\cong L_0$
and, by Fact \ref{T4100}(a), contains a copy of $L_0$. Thus, by (b), $f[L_0]\cap L_i^1\neq\emptyset$, for infinitely
many $i\in \o$. But this is impossible because $f[L_0]<f[L_1]$. Thus $f[L_0]\subset L_0$ and, hence, $f[L_0]\in \P (L_0)$.
By Proposition \ref{T4103}(a) we have $f[L_0]\cap L_i^0\neq\emptyset$, for infinitely
many $i\in \o$, which implies $f[L_1]\subset L_1$.
\kdok

\begin{prop}\rm\label{T4106}
(Finite sums of $\o$-sums) Let $L=\sum _{i\leq n } L_i$, where  $L_i\in \H$ are $\o$-sums
of sequences in $\H$ satisfying (\ref{EQ4100})
and $L_i +L_{i+1}\not\in \H$, for $i<n$. Then

(a) If $f:L\hookrightarrow L$, then $f[L_i]\subset L_i$, for each $i\leq n$;

(b) $\P (L)=\{ \bigcup _{i\leq n} C_i : \forall i\leq n \;\; C_i \in \P (L_i)\}$;

(c) sm$\langle \P (L) , \subset \rangle$ is a $\sigma$-closed pre-order.
\end{prop}

\dok
(a) For $n=1$ this is (c) of Lemma \ref{T4105}. Assuming that the statement is true for $n-1$ we prove that
it is true for $n$.
Suppose that $f[L_0]\not\subset L_0$. Then, since $f[L_n]\subset \bigcup _{i\leq n}L_i$, for
$i^* = \max \{ i\leq n: f[L_i]\not\subset \bigcup _{j\leq i}L_j \}$ we have $0\leq i^* <n$,
$f[L_{i^*}]\not\subset \bigcup _{j\leq i^*}L_j $ and $f[L_{i^*+1}]\subset \bigcup _{j\leq i^*}L_j \cup L_{i^* +1}$.
Since $f[L_{i^*}] < f[L_{i^*+1}]$ we have $f[L_{i^* +1}] \subset L_{i^*+1}$ so $f[L_{i^*}]\cap L_{i^*+1}$ is a final part
of $f[L_{i^*}]\cong L_{i^*}$ and, by Fact \ref{T4100}(a), contains a copy of $L_{i^*}$. This copy is contained in
the union of finitely many summands of $L_{i^* +1}$. But,
since $L_{i^*} +L_{i^*+1}\not\in \H$, this is impossible by Lemma \ref{T4105}(b).
Thus $f[L_0]\subset L_0$ and, by Proposition \ref{T4103}(a), the set $f[L_0]$ intersects infinitely many summands
of $L_0$, which implies $f[L_1 \cup \dots \cup L_n]\subset L_1 \cup \dots \cup L_n$. Thus, by the induction hypothesis,
$f[L_i]\subset L_i$, for each $i\in \{1, \dots , n\}$.

(b) The inclusion ``$\supset$" is evident and we prove ``$\subset$". If $C\in \P (L)$ and $f:L\hookrightarrow L$, where $C=f[L]$, then
by (a), $C_i =f[L_i]\subset L_i$ and, hence, $C_i \in \P (L_i)$ and, clearly, $C=\bigcup _{i\leq n}C_i$.

(c) By the statement (b) and, since the sets $L_i$, $i\leq n$, are disjoint, the mapping
$F:\prod _{i\leq n}\la \P (L_i), \subset \ra \rightarrow \la \P (L), \subset \ra$ given by
$F(\la C_0, \dots , C_n \ra )= C_0 \cup \dots \cup C_n$ is an isomorphism and, by Fact \ref{T2226},
$\sm \la \P (L),\subset \ra
\cong \sm (\prod _{i\leq n}\la \P (L_i), \subset \ra )
\cong \prod _{i\leq n}\sm \la \P (L_i), \subset \ra$.
By Proposition \ref{T4103}(c), the pre-orders $\sm \la \P (L_i), \subset \ra$, $i\leq n$, are $\sigma$-closed, and, by
Fact \ref{T4043} the same holds for their direct product and, hence, for sm$\langle \P (L) , \subset \rangle$ as well.
\kdok

\noindent
The following dual statements can be proved in the same way.

\begin{lem}\rm\label{T4107}
Let $L_0=\sum _{\o ^*} L_i^0 , L_1=\sum _{\o ^*} L_i^1 \in \H$, where $\la L_i^0 :i\in \o \ra$  and $\la L_i^1 :i\in \o \ra$
are sequences in $\H$ satisfying (\ref{EQ4100}). Then

(a) $\exists i\in \o \;\; L_1 \hookrightarrow L_i ^0 \Leftrightarrow
\exists m \in \o \;\; L_1 \hookrightarrow L_m^0 +\dots +L_0^0$;

(b) $L_0 +L_1 \not\in \H \Rightarrow \neg \exists m\in \o \;\;L_1 \hookrightarrow L_m^0 +\dots +L_0^0$.

(c) If $L=L_0 +L_1 \not\in \H $ and $f:L\hookrightarrow L$, then $f[L_k]\subset L_k$, for $k=0,1$.
\end{lem}
\begin{prop}\rm\label{T4108}
(Finite sums of $\o ^*$-sums) Let $L=\sum _{i<n } L_i$, where  $L_i\in \H$ are $\o ^*$-sums and
$L_i +L_{i+1}\not\in \H$, for $i<n-1$. Then

(a) If $f:L\hookrightarrow L$, then $f[L_i]\subset L_i$, for each $i<n$;

(b) $\P (L)=\{ \bigcup _{i<n} C_i : \forall i<n \;\; C_i \in \P (L_i)\}$;

(c) sm$\langle \P (L) , \subset \rangle$ is a $\sigma$-closed pre-order.
\end{prop}

\section{$\o ^*$-sum plus $\o $-sum}

\begin{lem}\rm\label{T4109}
Let $L=L_0 +L_1$, where $L_0=\sum _{\o ^*} L_i^0 , L_1=\sum _{\o } L_i^1 \in \H$ and
$\la L_i^0 :i\in \o \ra$  and $\la L_i^1 :i\in \o \ra$
are sequences in $\H$ satisfying (\ref{EQ4100}). Then

(a) $\exists i\in \o \;\; L_0 \hookrightarrow L_i ^1 \Leftrightarrow
\exists m \in \o \;\; L_0 \hookrightarrow L_0^1 + \dots + L_m^1$;

(b) $\exists i\in \o \;\; L_1 \hookrightarrow L_i ^0 \Leftrightarrow
\exists m \in \o \;\; L_1 \hookrightarrow L_m^0 + \dots + L_0^0$;

(c) If $L_0 +L_1 \not\in \H$, then
\begin{equation}\label{EQ4103}
\forall m\in \o \;\;(L_0 \not\hookrightarrow L_0^1 + \dots + L_m^1 \;\;\land \;\;
                     L_1 \not\hookrightarrow L_m^0 + \dots + L_0^0 ).
\end{equation}
\end{lem}
\dok
(a) If $f: L_0 \hookrightarrow \sum _{i\leq m}L_i ^1$ and
$i_0=\min \{i\leq m : f[L_0]\cap L_i^1\neq\emptyset\}$,
then $f[L_0]\cap L_{i_0}^1$ is a initial part of the ordering $f[L_0]\cong L_0$ and, by Fact \ref{T4100}(c),
contains a copy of $L_0$. Thus $L_0\hookrightarrow L_{i_0}^1$. The proof of (b) is dual.

(c) If $L_0 \hookrightarrow \sum _{i\leq m}L_i ^1$ then, by (a), there are $i_0\in \o$ and
$f: L_0\hookrightarrow L_{i_0}^1$. Then $\la L_0, L_0^1 , L_1^1 , \dots , L_{i_0}^1 ,\dots \ra$ is a sequence in $\H$
satisfying (\ref{EQ4100}) and, hence,  $L_0 + L_1=L_0 + L_0^1 + L_1^1 + \dots + L_{i_0}^1 +\dots \in \H$.
If $L_1 \hookrightarrow L_m^0 + \dots + L_0^0$, we prove $L_0 + L_1\in \H$ in a similar way.
\kdok

\begin{prop}\rm\label{T4110}
Let $L=L_0 +L_1 \not\in \H$, where $L_0=\sum _{\o ^*} L_i^0 , L_1=\sum _{\o } L_i^1 \in \H$ and
$\la L_i^0 :i\in \o \ra$  and $\la L_i^1 :i\in \o \ra$
are sequences in $\H$ satisfying (\ref{EQ4100}). Then

(a) $A\subset L$ contains a copy of $L$ iff for each $i, m\in \o $ there is a finite $K \subset \o \setminus m$ such that
$L_i^0 \hookrightarrow \bigcup _{j\in K }L_j^0 \cap A$ and $ L_i^1 \hookrightarrow \bigcup _{j\in K }L_j^1 \cap A$.

(b) If $A,B\in \P (L)$, then $A\leq B$  iff
for each $C\in \P (L)$ satisfying $C\subset A$ and each $i, m\in \o$ there is a finite $K \subset \o \setminus m$
such that
$L_i^0 \hookrightarrow \bigcup _{j\in K }L_j^0 \cap C \cap B$ and
$L_i^1 \hookrightarrow \bigcup _{j\in K }L_j^1 \cap C \cap B$.

(c) sm$\langle \P (L) , \subset \rangle$ is a $\sigma$-closed pre-order.
\end{prop}

\dok
(a) ($\Rightarrow$) Let $C\in \P (L)$, $C\subset A$, $f:L\hookrightarrow L$ and $C=f[L]$. First we prove
\begin{equation}\label{EQ4117}\textstyle
\exists C_0 \in \P (L_0)\;\; \exists C_1 \in \P (L_1)\;\; C_0\cup C_1 \subset A.
\end{equation}
Suppose that $f[L_0]\subset L_1$. Then, by Lemma \ref{T4109}(c), $f[L_0]\cap L_i ^1\neq \emptyset$, for infinitely
many $i\in \o$. But this is impossible since $f[L_0]<f[L_1]$. Thus $f[L_0]\cap L_0 \neq \emptyset$, this set is an
initial part of the order $f[L_0]\cong L_0$ and, by Fact \ref{T4100}(c), there is $C_0\in \P (L_0)$ such that
$C_0 \subset f[L_0]\cap L_0 \subset C\subset A$. Similarly, there is
$C_1\in \P (L_1)$ such that
$C_1 \subset f[L_1]\cap L_1 \subset C\subset A$ and (\ref{EQ4117}) is proved.

Let $i,m \in \o$.
By (\ref{EQ4117}) we have $C_0 \subset A\cap L_0 \subset L_0$ and $C_1 \subset A\cap L_1 \subset L_1$,
so, by Proposition \ref{T4103}(a), there are finite sets $K_0 ,K_1 \subset \o \setminus m$ such that
$L_i^0 \hookrightarrow \bigcup _{j\in K_0}L_j^0 \cap A \cap L_0$ and
$L_i^1 \hookrightarrow \bigcup _{j\in K_1}L_j^1 \cap A \cap L_1$. Clearly, $K=K_0\cup K_1$ is a finite subset
of $\o \setminus m$ and $L_i^0 \hookrightarrow \bigcup _{j\in K}L_j^0 \cap A $ and
$L_i^1 \hookrightarrow \bigcup _{j\in K}L_j^1 \cap A $.

($\Leftarrow$) Suppose that the given condition is satisfied by $A$. Then, by Proposition \ref{T4103}(a), there are
$C_0\in \P (L_0)$ and $C_1\in \P (L_1)$ such that $C_0\subset A\cap L_0$ and $C_1\subset A\cap L_1$. Now
$\P (L)\ni C_0 \cup C_1 \subset A$.

(b) By (\ref{EQ4141}), $A\leq B$ iff for each $C\in \P (L)$ satisfying $C\subset A$
the set $C\cap B$ contains a copy of $L$. Now we apply (a) to $C\cap B$.

(c)  For $A_n \in \P (L)$, $n\in \o$, where $A_0 \geq A_1 \geq \dots$ we will construct $A \in \P (L)$ such that
$A\leq A_n$, for all $n\in \o$.
First, by Fact \ref{T2226}(c), there are $C_i \in \P (L)$,
$i\in \o$, such that $C_0=A_0$ and
\begin{equation}\label{EQ4118}\textstyle
\forall i\in \o \;\; C_i \subset A_0 \cap \dots \cap A_i .
\end{equation}
By recursion we define
the sequences $\la K_i :i\in \o \ra $, $\la f_i^0 :i\in \o \ra $ and $\la f_i^1 :i\in \o \ra $
such that for each $i\in \o$

(i) $K_i \in [\o ]^{<\o }$,

(ii) $K_i < K_{i+1} $,

(iii) $f_i^0 : L_i^0 \hookrightarrow \bigcup _{j\in K_i}L_j^0 \cap C_i$,

(iv)  $f_i^1 : L_i^1 \hookrightarrow \bigcup _{j\in K_i}L_j^1 \cap C_i$.

\noindent
Since $C_0=A_0\in \P (L)$, by (a) (for $i=m=0$), there exist $K_0 \in [\o ]^{<\o }$,
$f_0^0 : L_0^0 \hookrightarrow \bigcup _{j\in K_0}L_j^0 \cap C_0$ and
$f_0^1 : L_0^1 \hookrightarrow \bigcup _{j\in K_0}L_j^1 \cap C_0$.

Let the sequences $K_0, \dots , K_{i'}$, $f_0^0, \dots , f_{i'}^0$ and $f_0^1, \dots , f_{i'}^1$ satisfy (i)-(iv).
Since $A_{i'+1}\leq A_{i'}$, $C_{i'+1}\in \P (L)$
and, by (\ref{EQ4118}), $C_{i'+1}\subset A_{i'+1}$, according to (b), for $i'+1$ and
$m=\max (K_0 \cup \dots \cup K_{i'})+1$ there are
\begin{equation}\label{EQ4119}\textstyle
K_{i'+1} \in [\o \setminus (\max (K_0 \cup \dots \cup K_{i'})+1 )]^{<\o }
\end{equation}
\begin{equation}\label{EQ4120}\textstyle
f_{i'+1}^0 : L_{i'+1}^0 \hookrightarrow \bigcup _{j\in K_{i'+1}}L_j^0 \cap C_{i'+1}
\end{equation}
\begin{equation}\label{EQ4121}\textstyle
f_{i'+1}^1 : L_{i'+1}^1 \hookrightarrow \bigcup _{j\in K_{i'+1}}L_j^1 \cap C_{i'+1}
\end{equation}
(since, by (\ref{EQ4118})), $C_{i'+1} \cap A_{i'}=C_{i'+1}$). By (\ref{EQ4119})) we have (i) and (ii).  (iii) and (iv)
follow from (\ref{EQ4120}) and (\ref{EQ4121}). The recursion works.

Let $f=\bigcup _{i\in \o }f_i^0 \cup \bigcup _{i\in \o }f_i^1$. By (ii) and (iii), $i_1 <i_2$ implies $K_{i_1}<K_{i_2}$,
which implies $f_{i_1}^0 [L_{i_1}^0]> f_{i_2}^0[L_{i_2}^0]$ and $f_{i_1}^1 [L_{i_1}^1]< f_{i_2}^1[L_{i_2}^1]$
and, hence, $f:L\hookrightarrow L$.
Thus
\begin{equation}\label{EQ4122}\textstyle
A=f[L]=\bigcup _{i\in \o }f_i^0[L_i^0] \cup \bigcup _{i\in \o }f_i^1 [L_i^1]\in \P (L).
\end{equation}
Using the characterization from (b), for $n^* \in \o$ we show that $A\leq A_{n^*}$.
So, for $C^*\in \P (L)$ such that $C^* \subset A$ and $i^*,m^* \in \o$ we prove that
\begin{equation}\label{EQ4123}\textstyle
\exists K \in [\o \setminus m^* ]^{<\o }\;(L_{i^*}^0\hookrightarrow \bigcup _{j\in K}L_j^0 \cap C^* \cap A_{n^*}
\land L_{i^*}^1\hookrightarrow \bigcup _{j\in K}L_j^1 \cap C^* \cap A_{n^*}).
\end{equation}
By (ii)-(iv) and (\ref{EQ4122}) we have $A=\sum _{\o ^*}\Lambda _i^0 + \sum _{\o }\Lambda _i^1 \cong L$,
where $\Lambda _i^0 =f_i^0[L_i^0]\cong L_i^0$ and $\Lambda _i^1 =f_i^1[L_i^1]\cong L_i^1$,
so $A$ is a sum of an $\o ^*$-sum, $\Lambda _0 =\sum _{\o ^*}\Lambda _i^0\cong L_0$
and an $\o $-sum, $\Lambda _1 =\sum _{\o }\Lambda _i^1 \cong L_1$. In addition, $L_0 + L_1\not\in \H$
implies $\Lambda_0 + \Lambda_1\not\in \H$.

Since $C^* \cong L \cong A$ and $C^* \subset A$ we have $C^* \in \P (A)$ so, applying (a) to the linear order
$A$ instead of $L$ we obtain
\begin{equation}\label{EQ4124}\textstyle
\forall i,m \in \o \;\; \exists K \in [\o \setminus m ]^{<\o }\;\;
(\Lambda _i^0 \hookrightarrow \bigcup _{j\in K}\Lambda _j^0 \cap C^*
\land \Lambda _i^1 \hookrightarrow \bigcup _{j\in K}\Lambda _j^1 \cap C^* ).
\end{equation}
Let $m' >m^* , n^*$. By (\ref{EQ4124}), for $i^*$ and $m'$ there is
\begin{equation}\label{EQ4125}\textstyle
K^* \in [\o \setminus m']^{<\o } \mbox{ such that}
\end{equation}
\begin{equation}\label{EQ4126}\textstyle
\Lambda _{i^*}^0 \hookrightarrow \bigcup _{j\in K^*}\Lambda _j^0 \cap C^*
\land \Lambda _{i^*}^1 \hookrightarrow \bigcup _{j\in K^*}\Lambda _j^1 \cap C^*
\end{equation}
By (\ref{EQ4125}), for $j\in K^*$ we have $j>n^*$ and, by (\ref{EQ4118}), $C_j \subset A_{n^*}$. Thus, by (iii) and (iv)
we have
$\Lambda _j^0 \subset \bigcup _{s\in K_j}L_s^0 \cap C_j \subset \bigcup _{s\in K_j}L_s^0 \cap A_{n^*}$
and
$\Lambda _j^1 \subset \bigcup _{s\in K_j}L_s^1 \cap C_j \subset \bigcup _{s\in K_j}L_s^1 \cap A_{n^*}$
which, together with
(iii),(iv) and (\ref{EQ4126}) gives
$L_{i^*}^0 \hookrightarrow \Lambda _{i^*}^0 \hookrightarrow  \bigcup _{j\in K^*}f_j[\Lambda _j^0] \cap C^*
\subset  \bigcup _{j\in K^*} \bigcup _{s\in K_j}L_s^0 \cap A_{n^*} \cap C^*
=\bigcup _{s\in \bigcup _{j\in K^*} K_j}L_s^0 \cap C^* \cap A_{n^*}$.
Similarly we prove that $L_{i^*}^0 \hookrightarrow \bigcup _{s\in \bigcup _{j\in K^*} K_j}L_s^0 \cap C^* \cap A_{n^*}$.

In order to finish the proof of (\ref{EQ4123}) we show that $\bigcup _{j\in K^*} K_j \cap m^* =\emptyset$.
By (\ref{EQ4125}), for $j\in K^*$ we have $j>m^*$. By (ii) the sequence $\la \min K_i : i\in \o \ra$ is increasing
and, hence, $\min K_j \geq j >m^*$, which implies $K_j \cap m^* =\emptyset$ and (\ref{EQ4123}) is proved.
\kdok

\section{The general case}

For $L\in \CS$, let $m(L)=\min \{ n\in \o : L \mbox{ is a sum of }n \mbox{ elements of } \H\}$. For $m\in \N$,
let $\CS _m =\{ L\in \CS : m(L)=m\}$.

\begin{lem}\rm\label{T4112}
(a) There is no $L\in \H$ such that $L=\sum _{\o ^*} L_i^0 $ and $L=\sum _{\o } L_i^1 $, where
$\la L_i^0 :i\in \o \ra$  and $\la L_i^1 :i\in \o \ra$
are sequences in $\H$ satisfying (\ref{EQ4100}).

(b) Let $L\in \CS _m$ and $L_0, \dots ,L_{m-1}\in \H$, where $L=L_0 +\dots +L_{m-1}$. Then
\begin{equation}\label{EQ4127}\textstyle
\forall i<m \;\; (|L_i|=1 \;\veebar\; L_i \mbox{ is an }\o \mbox{-sum} \;\veebar\; L_i \mbox{ is an }\o^* \mbox{-sum})
\end{equation}
\begin{equation}\label{EQ4128}\textstyle
\forall i<m-1 \;\; L_i +L_{i+1}\not\in \H .
\end{equation}
\begin{equation}\label{EQ4129}\textstyle
|L_i|=1 \Rightarrow (L_{i+1} \mbox{ is not an }\o \mbox{-sum } \;\; \land \;\;
L_{i-1} \mbox{ is not an }\o ^* \mbox{-sum }).
\end{equation}
\end{lem}
\dok
(a) On the contrary, by Fact \ref{T4100}, $L$ would be both left and right indecomposable and, for a partition
$L=L' +L''$ there would be $C',C''\cong L$ such that $C'\subset L'$ and $C'' \subset L''$, which would imply
$L+L\hookrightarrow L$. But this is impossible by Fact \ref{T4102}.

(b) The first statement follows from (a), the second from the minimality of $m$ and the third from the second statement
(1 + $\omega$-sum is an $\omega$-sum satisfying (\ref{EQ4100})).
\kdok

\begin{lem}\rm\label{T4111}
If $m\in \N$, $L\in \CS _m$, $L_0, \dots ,L_{m-1}\in \H$, where $L=L_0 +\dots +L_{m-1}$, and $f:L\hookrightarrow L$,
then for each $i<m$ there is $C_i \in \P (L_i)$ such that $C_i\subset f[L_i]$.
\end{lem}
\dok
We use induction. For $m=1$ the statement is trivially true.

Suppose that the statement holds for all $k\leq m$.
Let $L\in \CS _{m+1}$, $L_0, \dots ,L_m \in \H$, $L=L_0 + L_1 + \dots +L_m$ and $f:L\hookrightarrow L$.
Let $L'= L_1 + \dots +L_m$.

\vspace{2mm}
\noindent
{\it Claim 1}. $f[L_1]\cap L_0$ does not contain a copy of $L_1$.

\vspace{2mm}
\noindent
{\it Proof of Claim 1}.
On the contrary suppose that $L_1\cong C_1 \subset f[L_1]\cap L_0$.

First we show that $L_0$ is an $\o ^*$-sum.
Namely, $|L_0|=1$ would imply $C_1=L_0=f[L_1]$, which is impossible
because $f[L_0]< f[L_1]$. Suppose that $L_0$ is an $\o$-sum, $L_0=\sum _{\o }\Lambda _i$.
Then, since $f[L_0]< f[L_1]\cap L_0$, $L_0\hookrightarrow \sum _{i\leq m }\Lambda _i$, for some $m\in \o$,
which is impossible by Proposition \ref{T4103}(a).

Thus $L_0$ is an $\o ^*$-sum, $L_0=\sum _{\o ^* }L_i^0$  and, by (\ref{EQ4127}) and (\ref{EQ4129}),
$L_1$ is either an $\o$-sum or an $\o ^*$-sum.
Since $f[L_0]< f[L_1]\cap L_0 \hookleftarrow L_1$, there is $m\in \o $ such that
$L_1\hookrightarrow L_m^0 + \dots +L_0^0$. By (\ref{EQ4128}) we have $L_0+L_1\not\in \H$ and this is impossible
by Lemma \ref{T4109}(c) in the first case and Lemma \ref{T4107}(b) in the second. A contradiction.
Claim 1 is proved.

By (\ref{EQ4127}), regarding the summand $L_1$ we have the following three cases.

\vspace{2mm}
\noindent
{\it Case 1}: $|L_1|=1$.
Then, by Claim 1, $f[L_1]\cap L_0 =\emptyset$, which implies that
$f\upharpoonright L' : L'\hookrightarrow L'$. Clearly $m(L')\leq m$ and $m(L')< m$ is impossible, because of the
minimality of $m(L)$. Thus $m(L')= m$ and, by the induction hypothesis,
\begin{equation}\textstyle \label{EQ4130}
\forall i\in \{ 1,\dots m \}\;\; \exists C_i \in \P (L_i)\;\; C_i \subset (f\upharpoonright L')[L_i ]=f[L_i].
\end{equation}
Since $|L_1|=1$ we have
$C_1=L_1=f[L_1]>f[L_0]$, for $C_0=f[L_0]$ we have $C_0\in \P (L_0)$ and the proof is over.

\vspace{2mm}
\noindent
{\it Case 2}: $L_1$ is an $\o ^*$-sum.
By Fact \ref{T4100}(c), $f[L_1]\cap L_0 \neq\emptyset $ would imply that  $f[L_1]\cap L_0 $ contains a copy
of $L_1$, which is impossible by Claim 1. Thus $f[L_1]\cap L_0 = \emptyset $ and, as in Case 1, we have
(\ref{EQ4130}). In particular, $\P (L_1)\ni C_1 \subset f[L_1]$ and, by Proposition \ref{T4103}(a) (for $\o ^*$-sums),
$f[L_1]$ intersects infinitely many summands of $L_1$, which implies $f[L_0]\subset L_0$.
Again, for $C_0=f[L_0]$ we have $C_0\in \P (L_0)$ and the proof is over.

\vspace{2mm}
\noindent
{\it Case 3}: $L_1$ is an $\o $-sum. By (\ref{EQ4127}) and (\ref{EQ4129}),
regarding the summand $L_0$ we have the following two subcases.

\vspace{2mm}
\noindent
{\it Subcase 3.1}: $L_0$ is an $\o $-sum.
$f[L_1]\cap L_0 \neq\emptyset $ would imply that $L_0$ is embeddable in an initial part of $L_0$, which is impossible
by Proposition \ref{T4103}(a). Thus $f[L_1]\cap L_0 = \emptyset $ and, as in Case 1, we have (\ref{EQ4130}).
Since $C_1 \subset f[L_1]\cap L_1$ we have $f[L_0]\subset L_0 \cup L_1$. Suppose that $f[L_0]\cap L_1\neq\emptyset$.
Then $f[L_0]\cap L_1$ is contained in finitely many summands of $L_1$
and, by Fact \ref{T4100}(a), contains a copy of $L_0$, which is impossible by
(\ref{EQ4128}) and Lemma \ref{T4105}(b).
Thus $f[L_0]\subset L_0$ and, for $C_0=f[L_0]$ we have $C_0\in \P (L_0)$ which, together with (\ref{EQ4130}), finishes
the proof.

\vspace{2mm}
\noindent
{\it Subcase 3.2}: $L_0$ is an $\o ^*$-sum.
Let $L_0=\sum _{\o ^*}A_i$ and $L_1=\sum _{\o}B_i$. By Claim 1, there is $x\in L_1$ such that $L_0 <\{ f(x)\}$.
By Fact \ref{T4100}(b), there is $L_1'\cong L_1$ such that $L_1'\subset [x, \infty )_{L_1}$. Let
$\varphi : L_1 + L_2 + \dots + L_m \rightarrow L_1' + L_2 + \dots + L_m$ be an isomorphism, where
$\varphi \upharpoonright L_i = id _{L_i}$, for $i\in \{ 2,3, \dots , m \}$. Then $f\circ \varphi : L' \hookrightarrow L'$
and, by the induction hypothesis, there are $C_i \in \P (L_i)$, $i\in \{ 1,\dots m \}$,
satisfying $C_i \subset f[\varphi [L_i ]]$. Since $C_1 \subset f[\varphi [L_1]]=f[L_1']$ we have
\begin{equation}\textstyle \label{EQ4131}
C_1 \subset f[L_1' ]\cap L_1 \subset f[L_1]\cap L_1.
\end{equation}
\begin{equation}\textstyle \label{EQ4132}
\forall i\in \{ 2,\dots m \}\;\; (C_i \in \P (L_i) \; \land \; C_i \subset f[\varphi [L_i ]]=f[L_i]).
\end{equation}
By (\ref{EQ4131}) we have $f[L_0]\subset L_0 \cup L_1$.
Suppose that $f[L_0]\subset L_1$. Then, by (\ref{EQ4131}), $f[L_0]$ is contained in the union of finitely many
summands of $L_1$, which is impossible by (\ref{EQ4128}) and Lemma \ref{T4109}(c).
Thus $f[L_0]\cap L_0\neq \emptyset$ is an initial part of $f[L_0]\cong L_0$ and, by Fact \ref{T4100}(c), there is
$C_0 \cong L_0$ such that $C_0\subset f[L_0]\cap L_0$. By (\ref{EQ4131}) and (\ref{EQ4132}) the proof is over.
\kdok

Let $L\in \CS _m$ and $L_0, \dots ,L_{m-1}\in \H$, where $L=L_0 +\dots +L_{m-1}$. Then we have
(\ref{EQ4127}), (\ref{EQ4128}) and (\ref{EQ4129}) and
we divide $L$ into {\it blocks}, groups of consecutive summands $L_i$, in the following way:

- first we glue each two consecutive summands such that the
first is an $\o ^*$-sum and the second an $\o$-sum (blocks of the type D),

- then we divide the rest into the groups of consecutive (in $L$) $L_i$'s of the same form:
groups of singletons (blocks of the type A),
groups of $\o $-sums (blocks of the type B) and
groups of $\o ^*$-sums (blocks of the type C).

For example
$
111 | \o^*\o^* | \o^*\o |\o| 11 | \o^*\o |\o\o\o\o|\o^*\o^* .
$
More formally, we define a {\it block}  of $L$ as a sum of consecutive summands
$
B= L_i + L_{i+1} + \dots + L_{i+k} ,
$
where $k\geq 0$ and satisfying one of the following conditions.
\begin{itemize}
\item[(A)] $|L_j|=1$, for all $j\in \{ i, \dots , i+k \}$ and

         (i) $i=0 \lor |L_{i-1}|=\o $ and

         (ii) $i+k =m-1 \lor |L_{i+k+1}|=\o $;

\item[(B)] $L_j$ is an $\o $-sum, for all $j\in \{ i, \dots , i+k \}$ and

         (iii) $i=0 \lor (L_{i-1}$ is an $\o $-sum  $\land \;L_{i-2}$ is an $\o ^*$-sum) and

         (iv) $i+k=m-1 \; \lor \; L_{i+k+1}$ is not an $\o$-sum;

\item[(C)] $L_j$ is an $\o ^*$-sum, for all $j\in \{ i, \dots , i+k \}$ and

         (v) $i=0 \lor L_{i-1}$ is not an $\o ^*$-sum  and

         (vi) $i+k=m-1 \; \lor \; (L_{i+k+1}$ is an $\o ^*$-sum $\land \; L_{i+k+2}$ is an $\o $-sum);

\item[(D)] $k=1$ and $L_i$ is an $\o ^*$-sum and $L_{i+1}$ is an $\o $-sum.
\end{itemize}

\noindent
By $\Block (L)$ we will denote the set of blocks.

\begin{lem}\rm\label{T4113}
Blocks determine a partition of the set $\{ L_0, \dots , L_{m-1}\}$ and a partition of $L$ into convex parts.
\end{lem}
\dok
Let $L\in \CS _m$ and $L_0, \dots ,L_{m-1}\in \H$, where $L=L_0 +\dots +L_{m-1}$. First we show that
each summand $L_j$ is contained in some block. We have the following three cases

\vspace{2mm}
\noindent
{\it Case 1}: $|L_j|=1$.
Let $L_i, L_{i+1}, \dots , L_j, \dots L_{i+k}$ be the maximal sequence of consecutive summands
of size 1, including $L_j$. Then conditions (i) and (ii) are satisfied and, hence, $L_j$ belongs to a block of the type (A).

\vspace{2mm}
\noindent
{\it Case 2}: $L_j$ is an $\o$-sum.

\noindent
{\it Subcase 2.1}: $j=0$. Let $L_0, \dots ,L_k$ be a maximal sequence of consecutive $\o$-sums.
Then $k=m-1$ or $L_{k+1}$ is not an $\o$-sum so, conditions (iii) and (iv) are satisfied and $L_j$
belongs to a block of the type (B).

\noindent
{\it Subcase 2.2}: $j>0$ and $L_{j-1}$ is an $\o ^*$-sum. Then $L_{j-1}+L_j$ is a block of the type (D) containing $L_j$.

\noindent
{\it Subcase 2.3}: $j>0$ and $L_{j-1}$ is not an $\o ^*$-sum. Then, by (\ref{EQ4129}), $|L_{j-1}|\neq 1$, so,
by (\ref{EQ4127}), $L_{j-1}$ is an $\o$-sum. Let $L_i, L_{i+1}, \dots ,L_{j-1}, L_j, \dots ,L_{i+k} $
be the maximal sequence of consecutive $\o$-sums containing $L_j$. Then (iv) is true.

If $i=0$, then (iii) is true and $L_j$ belongs to a block of the type (B).

If $i>0$, then, by the maximality of the sequence and (\ref{EQ4129}) and (\ref{EQ4127}), $L_{i-1}$ is an $\o ^*$-sum.
Now $L_{i+1}, \dots ,L_{j-1}, L_j, \dots ,L_{i+k} $ satisfies (iii) and (iv), so it is a block of the type (B)
containing $L_j$ (since, clearly, $i+1\leq j$).

\vspace{2mm}
\noindent
{\it Case 3}: $L_j$ is an $\o ^*$-sum.

\noindent
{\it Subcase 3.1}: $j=m-1$. Let $L_i, \dots ,L_j$ be a maximal sequence of consecutive $\o ^*$-sums.
Then $i=0$ or $L_{i-1}$ is not an $\o ^*$-sum so, conditions (v) and (vi) are satisfied and $L_j$
belongs to a block of the type (C).

\noindent
{\it Subcase 3.2}: $j < m-1$ and $L_{j+1}$ is an $\o $-sum. Then $L_j+L_{j+1}$ is a block of the type (D) containing $L_j$.

\noindent
{\it Subcase 3.3}: $j < m-1$ and $L_{j+1}$ is not an $\o $-sum.
Since, by (\ref{EQ4129}), $|L_{j+1}|\neq 1$
by (\ref{EQ4127}) we have that $L_{j+1}$ is an $\o ^*$-sum.
Let $L_i, L_{i+1}, \dots ,L_j, L_{j+1}, \dots ,L_{i+k} $
be the maximal sequence of consecutive $\o ^*$-sums containing $L_j$.
Then (v) is true.

If $i+k=m-1$, then (vi) is true and $L_j$ belongs to a block of the type (C).

If $i+k<m-1$, then,
by the maximality of the sequence and (\ref{EQ4129}) and (\ref{EQ4127}), $L_{i+k+1}$ is an $\o $-sum.
Now $L_i, \dots ,L_{j-1}, L_j, \dots ,L_{i+k-1} $ satisfies (v) and (vi), so it is a block of the type (C)
containing $L_j$ (since, clearly, $j\leq i+k-1$).

Now we prove that different blocks are disjoint. Suppose that $B', B''\in \Block (L)$ and $x\in B'\cap B''$. Then
$x\in L_j$ for some $L_j$ contained in $B'\cap B''$. By (\ref{EQ4127}) we have the following three cases:

\vspace{2mm}
\noindent
{\it Case 1}: $|L_j|=1$. Then   $B'$ and $B''$ are blocks of the type (A). Since $L_j\subset B'\cap B''$,
by (i) and (ii) we have $B'=B''$.

\vspace{2mm}
\noindent
{\it Case 2}: $L_j$ is an $\o$-sum. Then, by Lemma \ref{T4112}(a), the blocks are of the type (B) or (D).

\noindent
{\it Subcase 2.1}: $B'$ and $B''$ are of the type (D). Then, since $L_j\subset B'\cap B''$ is an $\o$-sum, by
Lemma \ref{T4112}(a) we have $B'=B''$.

\noindent
{\it Subcase 2.2}: $B'$ and $B''$ are of the type (B). Then, since $L_j\subset B'\cap B''$, from (iii) and (iv)
it follows that in $L$ the blocks have the same beginning and the same end. Thus, $B'=B''$.

\noindent
{\it Subcase 2.2}: $B'$ is of the type (B) and $B''$ of the type (D). Then, by Lemma \ref{T4112}(a), $L_j$ is
the second summand of $B''$ and, hence, $B''=L_{j-1}+L_j$ and $B'=L_j + \dots +L_k$. But this is impossible by (iii)

\vspace{2mm}
\noindent
{\it Case 3}: $L_j$ is an $\o ^*$-sum. This case is dual to Case 2.
\kdok

\begin{lem}\rm\label{T4116}
If $m\in \N$, $L \in \CS _m$, $L_0, \dots ,L_{m-1}\in \H$, where $L=L_0 +\dots +L_{m-1}$
and $\Block (L)=\{ B_0, \dots B_r\}$,
then $\Block (L\setminus B_0)=\Block (L)\setminus \{ B_0\}$ .
\end{lem}
\dok
Let $L=L_0 +\dots +L_{n-1}+ L_n + \dots +L_{m-1}$, where $B_0 =L_0 +\dots +L_{n-1}$,
$L' =L\setminus B_0 =L_n + \dots +L_{m-1}$ and $0<n<m$. First we show that
\begin{equation}\textstyle \label{EQ4139}
\Block (L')\subset \Block(L).
\end{equation}
Let $B=L_i + \dots +L_{i+k} \in \Block (L')$. Clearly, if $B$ is of the type (D) in $L'$, then the same holds
in $L$ and $B\in \Block (L)$. If $B$ is of the type (A) (resp.\ (B), (C)), then it satisfies (ii) (resp.\ (iv), (vi)) in
$L'$ and, clearly, in $L$. If $i>n$, then, in addition, $B$ satisfies (i) (resp.\ (iii), (v)) in
$L'$ and, again, in $L$; thus $B\in \Block (L)$. So it remains to be proved that $B$ satisfies (i)
(resp.\ (iii), (v)) in $L$, when $i=n$.

\vspace{2mm}
\noindent
{\it Case 1}: $B$ is of the type (A). Then $|L_{n-1}|=1$ would imply that $B_0$ is not a block in $L$.
Thus $|L_{n-1}|=\o$ and $B$ satisfies (i) in $L$.

\vspace{2mm}
\noindent
{\it Case 2}: $B$ is of the type (B). Then $L_n$ is an $\o$-sum and, by (\ref{EQ4129}), $|L_{n-1}|=\o$.
By (iv) and (vi), $B_0$ is not of the type (B) or (C). Thus, $B_0$ is of the type (D) and, hence,  $B$
satisfies (iii) in $L$.

\vspace{2mm}
\noindent
{\it Case 3}: $B$ is of the type (C). Then $L_n$ is an $\o ^*$-sum.
Suppose that $L_{n-1}$ is an $\o ^*$-sum. Then $B_0$ must be of the type (C) and, by (vi) for $B_0$ in $L$, $L_{n+1}$
is an $\o$-sum. But then $B$ should be a block of the type (D) in $L'$, which is not true.
Thus $L_{n-1}$ is not an $\o ^*$-sum and, hence, $B$ satisfies (v) in $L$.

So (\ref{EQ4139}) is proved, which implies $\Block (L')\subset\Block (L)\setminus \{ B_0\} =\{ B_1, \dots B_r\}$.
By Lemma \ref{T4113} we have $\bigcup \Block (L')=L'=B_1 \cup \dots \cup B_r$, which gives the another inclusion.
\kdok
\begin{lem}\rm\label{T4114}
If $m\in \N$, $L \in \CS _m$, $L_0, \dots ,L_{m-1}\in \H$, where $L=L_0 +\dots +L_{m-1}$, and $f:L\hookrightarrow L$,
then for each $B\in$ Block$(L)$ we have $f[B]\subset B$.
\end{lem}

\dok
We prove the statement by induction. For $m=1$ it is trivially true.

Suppose that it is true for all $k< m$. Let $L=L_0 +\dots +L_{m-1}$ and $\Block (L)=\{ B_0, \dots B_r \}$.
If $r=0$, we are done. Otherwise we have
\begin{equation}\textstyle \label{EQ4133}
L= B_0 + L_{i+1} + \dots +L_{m-1} ,
\end{equation}
where $B_0= L_0 + \dots +L_i$. Let $L'=L_{i+1}+\dots +L_{m-1} $. By Lemma \ref{T4111},
\begin{equation}\textstyle \label{EQ4135}
\forall j\in \{ 0,\dots ,m-1 \} \;\; \exists C_j\in \P (L_j) \;\; C_j \subset f[L_j]\cap L_j .
\end{equation}
Regarding the type of $B_0$ we have the following cases.

\vspace{2mm}
\noindent
{\it Case 1}: $B_0$ is of the type (A).
Then, by (\ref{EQ4127}), (\ref{EQ4129}) and (ii), $L_{i+1}$ is an $\o ^*$-sum.
By (\ref{EQ4135}) and Proposition \ref{T4103}(a) (for $\o ^*$-sums),
$C_{i+1}$ intersects infinitely many summands of $L_{i+1}$ and, since $B_0$ is finite and $f[B_0]<f[L_{i+1}]$, we have
$f[B_0]=B_0$. Hence $f\upharpoonright L' : L' \hookrightarrow L' $ and $m(L')=m-i-1$. By Lemma \ref{T4116} we have
\begin{equation}\textstyle \label{EQ4134}
\Block (L')=\Block (L)\setminus \{ B_0\} =\{ B_1 , \dots ,B_r \}
\end{equation}
and, by the induction hypothesis, $f[B_j]= (f\upharpoonright L')[B_j]\subset B_j$, for $j>0$.

\vspace{2mm}
\noindent
{\it Case 2}: $B_0$ is of the type (B).
By Proposition \ref{T4103}(a) $C_i$ intersects infinitely many summands of $L_i$,
which implies that $f\upharpoonright L' : L' \hookrightarrow L' $.

If $|L_{i+1}|=1$, then $f[L_{i+1}]=L_{i+1}$ and, hence, $f[B_0]\subset B_0$. By (\ref{EQ4134}) and
the induction hypothesis $f[B_j]\subset B_j$, for $j>0$.

If $L_{i+1}$ is an $\o ^*$-sum, then $C_{i+1}$ intersects infinitely many summands of $L_{i+1}$ and, hence,
$f[B_0]\subset B_0$. Also, $C_i$ intersects infinitely many summands of $L_i$, which implies that
$f[L']\subset L'$. By (\ref{EQ4134}) and
the induction hypothesis $f[B_j]\subset B_j$, for $j>0$ again.

\vspace{2mm}
\noindent
{\it Case 3}: $B_0$ is of the type (C).
Then by (vi), $L_{i+1}$ is an $\o ^*$-sum.
By (\ref{EQ4135}) we have $C_{i+1}\subset f[L_{i+1}]\cap L_{i+1}$ and, by Proposition \ref{T4103},
$f[L_{i+1}]$ intersects infinitely many summands of $L_{i+1}$, which implies $f[B_0]\subset B_0$.
Suppose that $f[L_{i+1}]\cap L_i \neq\emptyset$.
By (\ref{EQ4135}), $C_i \subset f[L_i]\cap L_i$,
which implies that $f[L_{i+1}]\cap L_i$ is an initial part of $f[L_{i+1}]$
contained in an final part of $L_i$.
By Fact \ref{T4100}(c) $f[L_{i+1}]\cap L_i$ contains a copy of $L_{i+1}$,
which is impossible by Lemma \ref{T4107}(b) and (\ref{EQ4128}).
Thus $f[L_{i+1}]\cap L_i =\emptyset$, which implies $f[L']\subset L'$ and again, by (\ref{EQ4134}) and
the induction hypothesis $f[B_j]\subset B_j$, for $j>0$.

\vspace{2mm}
\noindent
{\it Case 4}: $B_0$ is of the type (D). Then $B_0=L_0+L_1$ and, by (\ref{EQ4135}) and Proposition \ref{T4103},
$f[L_1]$ intersects infinitely many summands of $L_1$, which implies
\begin{equation}\textstyle \label{EQ4136}
f[L']\subset L'.
\end{equation}
By (\ref{EQ4135}) there is $C_2$ such that
\begin{equation}\textstyle \label{EQ4137}
C_2 \in \P (L_2) \; \land \; C_2 \subset f[L_2]\cap L_2.
\end{equation}
Regarding the form of $L_2$ we distinguish the following three subcases.

$|L_2|=1$. Then, by (\ref{EQ4135}), $f[L_1]=L_1$ and, hence, $f[B_0]\subset B_0$ and we use (\ref{EQ4136}), (\ref{EQ4134})
and the induction hypothesis.

$L_2$ is an $\o ^*$-sum. By (\ref{EQ4137}) $f[L_2]$ intersects infinitely many summands of $L_2$ and, hence,
$f[B_0]\subset B_0$ and we use (\ref{EQ4136}), (\ref{EQ4134})
and the induction hypothesis.

$L_2$ is an $\o $-sum. By (\ref{EQ4137}) we have $f[L_1]\subset L_0 \cup L_1 \cup L_2$.
$f[L_1]\cap L_2 \neq\emptyset$ is impossible by Lemma \ref{T4105}(b), thus $f[B_0]\subset B_0$ and we continue as above.
\kdok

\begin{te}\rm\label{T4115}
For each $L\in \CS$, $\sm \la \P (L), \subset \ra$ is a $\sigma$-closed pre-order.
\end{te}
\dok
Let $L\in \CS _m$, $L=\sum _{i<r}B_i$, where $\Block (L)=\{ B_i : i<r\}$. First we prove
\begin{equation}\textstyle \label{EQ4138}
\P (L)=\{ \bigcup _{i<r } C_i : \forall i<r \;\; C_i \in \P (B_i )\}.
\end{equation}
The inclusion ``$\supset$" is evident. If $C\in \P (L)$, $f:L\hookrightarrow L$ and $C=f[L]$, then, by
Lemma \ref{T4114}, for $C_i=f[B_i]$, $i<r$, we  have $C_i \subset B_i$, $C_i \in \P (B_i)$ and $C=\bigcup _{i<r } C_i $
and ``$\subset$" holds as well.

Clearly, the mapping $F: \prod _{i<r}\la \P (B_i), \subset \ra \rightarrow \la \P (L), \subset \ra$ defined by
$$\textstyle
f(\la C_0 , \dots , C_{r-1} \ra) = \bigcup _{i<r } C_i
$$
is an isomorphism and, by Fact \ref{T2226}(d),(e)
$$\textstyle
\sm \la \P (L) , \subset \ra \cong \sm \prod _{i<r}\la \P (B_i), \subset \ra
= \prod _{i<r}\sm \la \P (B_i), \subset \ra.
$$
By Propositions \ref{T4106}, \ref{T4108} and \ref{T4110}
$\sm \la \P (B_i), \subset \ra$, $i<r$, are $\sigma$-closed partial orders and, by Fact \ref{T4043} their
product as well as the poset $\sm\la \P (L), \subset \ra$ is $\sigma$-closed.
\kdok

\section{Forcing by copies of countable scattered linear orders}
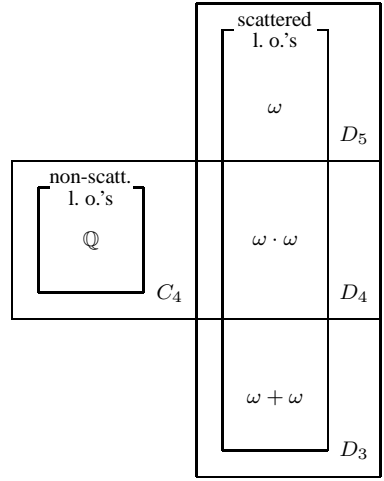
\begin{figure}[htb]
\begin{center}
\unitlength 0.7mm%1mm % = .854pt
\linethickness{0.4pt}
\ifx\plotpoint\undefined\newsavebox{\plotpoint}\fi % GNUPLOT compatibility
\begin{picture}(75,100)(0,0)
%----------------------------- linije --------------------------
\put(35,5){\line(1,0){35}}%1
\put(40,10){\line(1,0){20}}%2
\put(0,35){\line(1,0){70}}%3
\put(5,40){\line(1,0){20}}%4
\put(0,65){\line(1,0){70}}%5
\put(35,95){\line(1,0){35}}%6
\put(0,35){\line(0,1){30}}%7
\put(5,40){\line(0,1){20}}%8
\put(25,40){\line(0,1){20}}%9
\put(35,5){\line(0,1){90}}%10
\put(40,10){\line(0,1){80}}%11
\put(60,10){\line(0,1){80}}%12
\put(70,5){\line(0,1){90}}%13
\put(5,60){\line(1,0){2}}%14a
\put(23,60){\line(1,0){2}}%14b
\put(40,90){\line(1,0){2}}%15a
\put(58,90){\line(1,0){2}}%15b
%----------------------------- tacke --------------------------
%\put(70,10){\circle*{1}}%1
%----------------------------- tekst -----------------------
%\small
%\footnotesize
\scriptsize
%\tiny
\put(50,92){\makebox(0,0)[cc]{scattered}}%
\put(50,88){\makebox(0,0)[cc]{l.\ o.'s}}%
\put(15,62){\makebox(0,0)[cc]{non-scatt.}}%
\put(15,58){\makebox(0,0)[cc]{l.\ o.'s}}%
\put(15,50){\makebox(0,0)[cc]{$\Q$}}%
\put(50,75){\makebox(0,0)[cc]{$\o$}}%
\put(50,50){\makebox(0,0)[cc]{$\o \cdot \o$}}%
\put(50,20){\makebox(0,0)[cc]{$\o + \o$}}%
%--------------------A1D5------------------------------------------------
\put(30,40){\makebox(0,0)[cc]{$C_4$}}%
\put(65,10){\makebox(0,0)[cc]{$D_3$}}%
\put(65,40){\makebox(0,0)[cc]{$D_4$}}%
\put(65,70){\makebox(0,0)[cc]{$D_5$}}%
\end{picture}
\end{center}

\vspace{-7mm}
\caption{Countable linear orders}\label{F4002}
\end{figure}

The position of countable linear orders in Diagram \ref{F4001} is presented in Diagram \ref{F4002}.
%=============================================================
%===================================================================

By Theorem \ref{T4121} and Fact \ref{T4056}, CH implies that all posets of the form $\la \P (L), \subset \ra$, where $L$ is a scattered countable linear order,
are forcing equivalent to $(P(\o )/\Fin )^+$.
The following examples show that this is not true in general and that
the result of Theorem \ref{T4121} is the best possible: ``$\s$-closed" can not be replaced by ``$\o _2$-closed".
\begin{ex} \rm \label{EX4000}
It is consistent that the poset $\la \P (\o + \o ), \subset \ra$ is not
${\mathfrak h}$-distributive and, hence, not forcing equivalent to $(P(\o )/\Fin )^+$.

By Proposition \ref{T4106}, for $L = \o + \o$ the partial order $\la \P (L), \subset \ra$
is isomorphic to the product $\la [\o ]^\o, \subset \ra \times \la [\o ]^\o, \subset \ra$ and, by Fact \ref{T4091}(a),
$
\sq \la \P (\o + \o ), \subset \ra \cong (P(\o )/\Fin )^+ \times (P(\o )/\Fin )^+.
$
Now, by the result of Shelah and Spinas (Fact \ref{T4091}(b)), we have Con(${\mathfrak h}_2  < {\mathfrak h}$).
\end{ex}
\begin{ex} \rm \label{EX4001}
The poset $\sq \la \P (\o \cdot \o ), \subset \ra$ is not $\o _2$-closed and it is consistent that $\sq \la \P (\o \cdot \o ), \subset \ra$ is not
${\mathfrak h}$-distributive.
Clearly $\o \cdot \o \cong \la L, < \ra$, where
$L= \o \times \o $ and $\la i_0 , j_0 \ra < \la i_1 , j_1 \ra \Leftrightarrow i_0 < i_1 \lor (i_0 =i_1 \land j_0 < j_1)$.
Now $L=\sum _{i\in \o}L_i$, where $L_i=\{ i\} \times \o$
and first we show that $\P (L)=(\Fin \times \Fin )^+$.
By Proposition \ref{T4103}(a), if $A\in \P (L)$, then for each $m\in \o$ there is a finite set $K \subset \o \setminus m$ such that
$\o \hookrightarrow \bigcup _{i\in K}A\cap L_i$ and, hence, there is $i\geq m$  satisfying $|A \cap L_i | =\o$. Thus $A\not\in \Fin \times \Fin$.
Conversely, if $A\not\in \Fin \times \Fin$ and $\{ i\in \o : |A \cap L_i | =\o \} = \{ n_j : j\in \o\}$, where $n_0 < n_1 < \dots $,
then $A = \bigcup _{j\in \o}\Lambda _j$, where
$\Lambda _0 = \bigcup _{i\leq n_0} (A\cap L_i)$ and $\Lambda _j = \bigcup _{n_{j-1}< i\leq n_j} (A\cap L_i)$, for $j>0$. Clearly we have
$\Lambda _j \cong \o$ and, hence, $A\in \P (L)$.
So, $\la \P (L), \subset \ra = \la (\Fin \times \Fin )^+ ,\subset \ra$ and, by Fact \ref{T2226}(f),
$
\sq \la \P (\o \cdot \o), \subset \ra \cong  (P(\o \times \o)/(\Fin \times \Fin ))^+   .
$
Now we apply the results of Szyma\'nski and Zhou  and of Hern\'andez-Hern\'andez (Fact \ref{T4091}(c) and (d)).
\end{ex}
Some forcing-related properties of the posets $\sq \la \P (L ), \subset \ra$ are described in the following table.

\begin{center}
{%\small
\footnotesize
%\scriptsize
\begin{tabular}{c|c|c|c}
$L$          & $\sq \langle  {\mathbb P}(L), \subset \rangle $ is     & $\mathop{\rm sq}\nolimits  \langle  {\mathbb P} (L), \subset \rangle $ is & ZFC $\vdash \mathop{\rm sq}\nolimits  \langle  {\mathbb P} (L), \subset \rangle $   \\%[2mm]
                          &  isomorphic to             &                                    & is ${\mathfrak h}$-distributive              \\[2mm] \hline
                          &                                    &                                    &                                              \\%[2mm]
$\omega $             & $(P(\omega )/\mathop{\rm Fin}\nolimits )^+$                                                    & ${\mathfrak t}$-closed   & yes \\[2mm]
$\omega + \omega $    & $(P(\omega )/\mathop{\rm Fin}\nolimits )^+  \times  (P(\omega )/\mathop{\rm Fin}\nolimits )^+$ & ${\mathfrak t}$-closed   & no \\[2mm]
$\omega \cdot \omega$ & $(P(\omega \times \omega)/(\mathop{\rm Fin}\nolimits  \times \mathop{\rm Fin}\nolimits  ))^+$  & $\omega _1$ but not $\omega _2$-closed & no
\end{tabular}
}
\end{center}

\begin{rem}\rm \label{R4001}
Concerning Theorem \ref{T4121} we note that for countable ordinals we have more information. Namely,
by \cite{Kord}, if  $\a=\o ^{\g _n +r_n }s_n + \dots + \o ^{ \g _0 +r_0 }s_0 + k $ is a countable ordinal presented in the Cantor normal form, where $k\in \o$, $r_i \in \o$, $s_i \in \N$, $\g _i \in \Lim \cup \{ 1 \}$ and
$\g _n +r_n > \dots > \g _0 +r_0$,
then
\begin{equation}\label{EQ4078}\textstyle
\sq \la \P (\a ), \subset \ra \cong \prod _{i=0}^n \Big( \Big( \rp ^{r_i}\big( P(\o ^{\g _i} )/ \I _{\o ^{\g _i} }\big)\Big)^+ \Big)^{s_i},
\end{equation}
where, for an ordinal $\b$,
%$\I _\b$ denotes the collection of subsets of $\b$ which do not contain a copy of $\b$
$\I _\b =\{ C\subset \b : \b \not\hookrightarrow C\} $
and, for a poset $\P$, $\rp (\P )$ denotes the reduced
power $\P ^\o / \equiv _{\Fin }$ and $\rp ^{k+1}(\P )= \rp (\rp ^k (\P ))$.
In particular, for $\o \leq \a <\o ^\o$ we have
\begin{equation}\label{EQ4079}\textstyle
\sq \Big( \P \big(\sum _{i=n}^0 \o ^{1+r_i} s_i \big), \subset \Big)
\cong\prod _{i=0}^n \Big(\Big( \rp ^{r_i}\big( P(\o )/ \Fin \big)\Big)^+ \Big)^{s_i}.
\end{equation}
\end{rem}
\begin{rem}\rm \label{R4000}
By \cite{Kemb}, all countable equivalence relations, disconnected ultrahomogeneous graphs and disjoint unions of ordinals $\leq \o$ are in column $D$ of Diagram
\ref{F4001} as well. In addition, the corresponding posets of copies are  forcing equivalent to one of the following posets:

$((P(\omega )/\mathop{\rm Fin}\nolimits  )^+)^n$,  for some $n\in\N$,

$(P(\omega \times \omega)/(\mathop{\rm Fin}\nolimits  \times \mathop{\rm Fin}\nolimits  ))^+$, 

$(P(\Delta )/{\mathcal E}{\mathcal D}_{\mathrm{fin}}  )^+  \times  ((P(\omega )/\mathop{\rm Fin}\nolimits )^+)^n$, for some $n\in\omega$, 

\noindent
where $\Delta =\{ \langle  m,n\rangle \in {\mathbb N} \times {\mathbb N} : n\leq m\}$ and the ideal
${\mathcal E}{\mathcal D}_{\mathrm{fin}} \subset P(\Delta)$ is defined by:
$${\mathcal E}{\mathcal D}_{\mathrm{fin}}  = \{ S\subset \Delta :
\exists r\in {\mathbb N} \;\; \forall m \in {\mathbb N} \;\; |S\cap (\{ m \} \times \{ 1,2,\dots ,m\})|\leq r  \}.$$
\end{rem}

\end{document}